\documentclass[review]{elsarticle}
\usepackage{amsmath}
\usepackage{amssymb}  
\usepackage{mathtools,bbm}
\usepackage{color}
\usepackage[hyphens,spaces,obeyspaces]{url}
\usepackage{caption}

\newcommand{\bld}[1]{\boldsymbol{#1}}
\newcommand{\trn}[1]{{#1}^{\textrm{T}}}
\newcommand{\Ref}[1]{(\ref{#1})}
\newcommand{\rn}[1]{\mathbb{R}^{#1}_{+}}
\newcommand{\rng}[1]{\mathbb{R}^{#1}}
\newcommand{\Def}[1]{Definition \ref{#1}}
\newcommand{\cset}[2]{{\textup{Conset}}_{#1}(#2)}

\newcommand{\poly}[2]{\ensuremath{#1^{#2}-\alpha_{#2-1}#1^{#2-1}-\dots-\alpha_{1}#1-\alpha_{0}}}

\newtheorem{thm}{Theorem}
\newtheorem{cor}{Corollary}

\newtheorem{lem}{Lemma}
\newdefinition{prop}{Proposition}
\newdefinition{rmk}{Remark}
\newdefinition{defn}{Definition}
\newdefinition{exmp}{Example}
\newproof{pf}{Proof}
\newproof{pot}{Proof of Theorem \ref{thm2}}

\begin{document}

\begin{frontmatter}

\title{Controllability of Linear Positive Systems: An Alternative Formulation\tnoteref{mytitlenote}}
\tnotetext[mytitlenote]{Full draft for review only.}

\author[mymainaddress]{Yashar Zeinaly\corref{DCSC}}
\cortext[DCSC]{Delft Center for Systems and Control, Delft University of Technology}
\ead{y.zeinaly@tudelft.nl}

\author[mysecondaryaddress]{Jan H. van Schuppen\corref{DIAM}}
\cortext[DIAM]{Delft Institute of Applied Mathematics, Delft University of Technology}
\ead{J.H.vanSchuppen@tudelft.nl}

\author[mymainaddress]{Bart De Schutter\corref{DCSC}}
\ead{b.deschutter@tudelft.nl}

\address[mymainaddress]{Mekelweg 2, 2628 CD Delft, The Netherlands}
\address[mysecondaryaddress]{Mekelweg 4, 2628 CD Delft, The Netherlands}

\begin{abstract}
An alternative formulation for the controllability problem of single input linear positive systems is presented. Driven by many industrial applications, this formulations focuses on the case where the region of interest is only a subset of positive orthant rather than the entire positive orthant. To this end, we discuss the geometry of controllable subsets and develop numerically verifiable conditions for polyhedrality of controllable subsets. Finally, we provide a method to check for controllability of a target set based on our approach.
\end{abstract}
\begin{keyword}
Linear positive systems \sep Controllability \sep Cones \sep Polyhedral cones
\end{keyword}
\end{frontmatter}
\section{Introduction}\label{sect:intoduction}
In this paper, we revisit the ``controllability" concept for discrete-time linear positive systems. Motivated by applications underlying positive systems, we will re-define this concept. We will then provide necessary and sufficient conditions for controllability of a certain class of discrete-time linear positive systems.

The concept of positive systems arises in many applications such as econometrics \cite{Economics_Book}, bio-chemical reactors \cite{Compartmental_Book,Bioreactor}, compartmental systems \cite{vandenhof1998,haddad:chellaboina:hui}, and transportation system \cite{Zeinaly2014,Traffic_Network}, to name a few. The variables in such systems represent growth rates, concentration levels, mass accumulation, or flows, etc. Obviously, variables of this nature can only assume non-negative values.  The theory of positive dynamical systems has been developed to deal with this sort of systems. Of particular interest is the theory of linear positive systems \cite{Berman_book}, which has its roots in the theory of non-negative matrices and in the geometry of cones \cite{Plemmons_book,Davis1954,Gale,Vandergraft1968}. While the theory of linear positive systems has overlaps with general theory of linear systems, there are distinct differences between the two. This is due to the fact that linear positive systems are defined over a cone rather than over a linear subspace. Therefore, many properties of linear systems cannot be generalized to linear positive systems without proper treatment. Moreover, some concepts of general linear systems theory might have to be redefined for linear positive systems. One such property is the notion of ``controllability'' for linear positive systems. 

In many industrial applications one might be interested in investigating whether a certain state (e.g., concentration levels) can be reached by applying an appropriate control input. More generally one might be interested in characterizing all states that can be reached from a given initial state using nonnegative control inputs. With respect to this point of view, the alternative approach in this paper is based on the following key problem: \emph{Given a set of states, possibly a singleton, in $\rn{n}$, can the system initially at rest be steered in finite time to any state of the considered set by applying nonnegative control signal?}

The controllability of discrete-time linear positive systems has been widely studied in the literature. In most of the literature, it has been emphasized that the characterization of controllability for discrete-time linear positive systems takes a very peculiar form, which is very different from its counterpart for discrete-time linear systems \cite{Farina_book,Coxon_87,Bru2000}. Unlike discrete-time linear systems in $\rng{n}$ for which reachability is equivalent to reachability in $n$ steps \cite{Luenberger_book}, for discrete-time linear positive systems this does not hold and the timing issue becomes very critical, as noted in \cite{Coxon_87}, where they illustrate this using the model of  a pharmacokinetic system. However, inspired by the definition of reachability within the context of linear systems, most papers in the literature investigate and discuss necessary and sufficient conditions under which the positive orthant $\rn{n}$ is reachable. Among others, \cite{Bru2000,Valcher_96,Fanti1990,Guiver2014,Bartosiewicz:2013}, are some of the significant works that fall in this category. In \cite{Valcher_96,Bru2000} controllability of discrete-time linear positive systems is characterized using a graph-theoretic approach, and canonical controllability forms are derived as well. The authors of \cite{Guiver2014} have established a link between positive state controllability and positive input controllability of a related system, which is then used to obtain a controllability criterion. A good survey of similar results is provided in \cite{Rumchev_2000,Kaczorek2004}. Controllability results for special classes of 1D and 2D systems are provided in \cite{Positive2D_book}.

In this paper we first define and characterize the controllable subsets. Then, in Proposition~\ref{prop:polyhedral_conset_infty} and Proposition~\ref{prop:polyhedral_conset_f}, we present necessary and sufficient conditions for polyhedrality of the controllable subsets. Theorem~\ref{thm:polyherdal_conset_infty_characterization} and Theorem~\ref{thm:polyherdal_conset_f_characterization} provide a numerically verifiable method to check for polyhedrdality of the controllable subsets based on spectrum of $\bld{A}$. Finally, in Proposition~\ref{prop:LP_Solution} we propose a method to check for controllability of a given subset of $\rn{n}$.
The rest of this paper is organized as follows. In Section~\ref{sect:problem_Formulation}, inspired by the aforementioned application domains, we formally introduce our view of the controllability problem. In Section \ref{sect:Notation}, we introduce some notation that will be used in the sequel. A characterization of controllable subsets is then provided in Section~\ref{sect:characterization_Controllable_Subsets}, and the controllability problem is characterized in Section~\ref{sect:characterization_Controllability}.
\section{Problem Formulation}\label{sect:problem_Formulation}
\subsection{Classical view}
We will now introduce the classical view of the controllability problem as discussed in the literature highlighting that the stated conditions for controllability are often too strict and impractical. Then we will formally introduce our view of the controllability problem arguing why it is more suitable, especially from the application point of view.
\begin{rmk}
Different terminologies have been used for the concept of controllability of linear positive systems in the literature. Investigating whether a state is reachable from the origin has been referred to both as ``reachability'' and ``controllability from the origin." In this paper, in line with the latter terminology, since we assume the system is initially at rest, we will use controllability to refer to ``controllability from the origin."
\end{rmk} 
In most papers of the literature, the characterization of controllability of linear positive systems is based on the following definition, see \cite{Farina_book}.
\begin{defn}\label{def:Farina}
``A positive system is said to be completely reachable if all states $x\geq 0$ are reachable in finite time from the origin, that is, if $X_r=\mathbb{R}_+^n$,'' where $X_r=\mathbb{R}_+^n$ denotes the cone of all reachable states in finite time using nonnegative inputs.
\end{defn}
The underlying idea behind \Def{def:Farina} probably originates from making an analogy to reachability of linear systems. This definition is based on the assumption that the state space is $X=\rn{n}$. But if the system starts at the zero state, then it may not be possible with the existing inputs to reach all states of the system. Therefore the states to be reached may be restricted from the full positive orthant $X=\rn{n}$ to a smaller subset of the positive orthant. Hence the condition of controllability has to be adjusted as described in the remainder of the paper. The following theorem (\cite[Th. 27]{Farina_book}), states the necessary and sufficient condition for reachability with respect to \Def{def:Farina} for the single-input case.
\begin{thm}\label{thrm:Farina}
``A discrete-time positive system is completely reachable if it is possible to reorder its state variables in such a way that the input $u$ directly influences only $x_1$, and $x_i$ directly influences $x_{i+1}$ for $i=1,2,\dots,n-1$.''
\end{thm}
The results for the multi-input case based on \Def{def:Farina} are more involved, but they require that the matrix $[B,~AB,~\ldots,~A^{k}B]$ includes a monomial submatrix of dimension $n$, for some $k\in \mathbb{N}_{+}$ \cite{Valcher_96,Rumchev_2000,Fanti1990,Bru2000}. Such conditions are often too strong to be satisfied by most of practical systems. In addition, especially from the application point of view, complete reachability according to \Def{def:Farina} is not required in most of the cases since many practical positive systems operate in a constrained space, which is a strict subset of $\rn{n}$ and/or we are only interested in reachabiliy of states within a constrained space. For example in economical systems, one would be interested to know whether a certain growth rate can be achieved, which corresponds to checking whether a certain extremal ray is reachable. In bio-chemical reactors, it might be of interest to know whether a set of desired mass concentrations can be achieved by manipulating the inputs (e.g., flow of material). The set of desired concentrations is normally a small subset of $\mathbb{R}^n_+$.
\begin{exmp}\label{ex:0} Consider the discrete-time time-invariant linear positive system \begin{align}
\bld{x}(t+1)= &\ A\bld{x}(t)+B\bld{u}(t),\ \bld{x}(0)=\bld{x}_0,
\end{align}
with
$$A=\left[
\begin{matrix}
4 & 4\\11 & 2
\end{matrix}
\right],\
b = \left[\begin{matrix}
2\\1
\end{matrix}\right],\ \bld{x}_0=0.$$
It is of interest to determine whether the states in the cone $K\subset \rn{2}$, defined by \Ref{eq:example_1_1} and illustrated by Fig.~\ref{fig:E0}, can be reached in finite time:
\begin{equation}\label{eq:example_1_1}
K:\left\{\begin{aligned}
&3x_1-2x_2\geq 0,\\
&3x_2-2x_1\geq 0,\\
&x_1\geq 0,\ x_2\geq 0.\\ 
\end{aligned}\right.
\end{equation}
\begin{figure}
\centering
\includegraphics[width=0.3\linewidth]{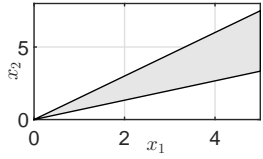}
\caption{\label{fig:E0}Example \ref{ex:0}. The shaded area, associated with $K$, represents the region of interest for which controllability needs to be checked.}
\end{figure}
Since $K\subset \rn{2}$, in order to answer this question using the classical approach, one needs to check the reachability of $\rn{2}$, which is a very conservative considering the fact that $K$ ``occupies'' only a small portion of $\rn{2}$. It can be verified that 
$$[b,Ab,\dots,A^kb]=\left[\begin{matrix}
2 & 12 & \cdots\\
1 & 24 & \cdots
\end{matrix}\right]$$ does not include a monomial submatrix of dimension $2$ for any $k\in \mathbb{N}_{+}$. Therefore, the conditions of Theorem \ref{thrm:Farina} do not hold and we cannot deduce anything about the reachability of $K$. Nevertheless, it will be later shown that $K$ is reachable from the origin.
\end{exmp}
\subsection{The approach of this paper}
From a practical point of view, the controllability problem boils down to whether it is possible to steer the system at rest to a given target set in finite time; and if this is the case, how long will it take to drive the system there. The controllable subset is defined as the subset of the state set containing those states that are reachable by either a finite or an infinite length nonnegative input signal. That subset of the state set is then a cone. Controllability is then defined as the requirement that the controllable subset contains the target set, which could be different from the positive orthant itself. Therefore the view point has to be changed by focusing on the controllable set, characterizing it, and the determination of conditions which guarantee that a particular subset of the positive orthant is contained in the controllable set. In addition, it will be shown that the controllable set in general does not have a finite characterization.\\
The notation of a linear positive system is formally defined in Section~\ref{sect:Notation}. In this section, the controllable subset is denoted as $\text{Conset}_{k}(x_{0})$, $\text{Conset}_{\textup{f}}(x_{0})$, or $\text{Conset}_{\infty}(x_{0})$ depending on whether the input sequence contains $k\in \mathbb{N}$ elements, a finite number, or an infinite number of elements.\\
Consider a discrete time-invariant linear positive system. The controllability problem is them composed of the following subproblems:
\begin{enumerate}[1.]
\item Characterize the controllable subsets $\text{Conset}_k(\bld{x}_0)$, $\text{Conset}_{\textup{f}}(\bld{x}_0)$, and\\ $\text{Conset}_{\infty}(x_{0})$ for the initial state $x_0=0$.
\item Determine whether or not the controllable sets $\text{Conset}_{\textup{f}}(x_{0})$ and $\text{Conset}_{\infty}(x_{0})$  can be computed in a finite number of steps.
\item Determine conditions on the system such that $\text{Conset}_{\infty}(x_{0})=\rn{n}$.
\item Considering a cone $C_{\text{obj}}\subseteq \rn{n}$ of control objectives or a subset of $\rn{n}$, determine sufficient and necessary conditions with respect to which the following condition holds: $C_{\text{obj}}\subseteq \text{Conset}_{\textup{f}}(0)$.
\end{enumerate}
\section{Concepts of Linear Positive Systems}\label{sect:Notation}
\subsection{Positive Real Numbers and Positive Matrices}\label{subsect:Positive_Marices}
The reader is informed of the following books on positive real numbers and positive matrices: \cite{Plemmons_book,Minc_book}. Books on positive systems or books with chapters on positive systems include \cite{Berman_book,Farina_book,Luenberger_book,Positive2D_book}.
The reader is assumed to be familiar with the integers, the real numbers, and vector spaces. Denote the set of the integers by $\mathbb{Z}$, the strictly positive integers by $\mathbb{Z}_{+} =\{1,2,\dots\}$, and the set of the natural numbers by $\mathbb{N} = \{0,1,2,\dots\}$. For any $n\in \mathbb{Z}_{+}$ denote the set of the first $n$ integers and of the first $n$ natural numbers by, respectively, $\mathbb{Z}_{n}=\{1,2,\dots,n\}$ and $\mathbb{N}_n =\{0,1,\dots,n\}$.

The real numbers are denoted by $\mathbb{R}$, the set of the positive real numbers by $\mathbb{R}_{+}=[0,\infty)$, and the set of the strictly positive real numbers by $\mathbb{R}_{s+}=(0,\infty)$. The $n$-dimensional vector space of tuples of real numbers is denoted by $\mathbb{R}^{n}$. The associated field of scalars is the set of the real numbers.

The set of the positive real numbers is a semi-ring. It is closed with respect to addition and with respect to multiplication. But it is not closed with respect to the inverse of addition (subtraction). The set of the strictly positive real numbers is closed with respect to inversion.

Consider the set of $n$ tuples of the positive real numbers $\rn{n}$, with the set of the positive numbers as the set of scalars. This set is closed with respect to addition but it does not have an inverse with respect to addition. The algebraic structure of $(\mathbb{R}_{+},\rn{n})$ is a semi-ring.

For a finite subset $S\subseteq \mathbb{R}^n$, $K\subseteq \mathbb{R}^n$ is the polyhedral cone generated by $S$ if it consists of all finite nonnegative linear combinations of elements of $S$. For a matrix $\bld{M}\in \mathbb{R}^{n\times m}$, we denote $\textup{cone}(\bld{M})$ as the cone generated by columns of $\bld{M}$. A {\em ray} of a cone is a line starting in the vertex of the cone and extending
to infinity, and lying on the boundary of the cone. It is called an {\em extreme ray} if it cannot be written as the convex combination of two other rays. A {\em polyhedral cone} is a cone for which there exists a finite number of extreme rays such that any vector starting at the vertex of the cone and extending to infinity, is a finite nonnegative linear combination of the extremal rays. A cone which is not polyhedral is also called a {\em round cone}. Thus, a cone is a round cone if there exists a non-denumerable number of extreme rays. An example of a round cone is the well known ice cream cone which may be found 
in \cite{Plemmons_book}.

For a finite set of complex numbers $S=\{s_1,s_2,\dots,s_k\}$, we denote $\rho(S)=\underset{s\in S}{\textup{max}}|s|$. For $\bld{A}\in \mathbb{R}^{n\times n}$, $\rho(\bld{A})=\rho({\textup{spec}(\bld{A})})$ is the spectral radius of $\bld{A}$, where $\textrm{spec}(\bld{A})$ denotes the set of its eigenvalues. We define the dominant subset of $S$ as $\sigma^{\rho}(S)=\{s\in \mathbb{C},|s|=\rho(S)\}$, and the non-dominant subset as $\sigma^{-}(S)=\{s\in \mathbb{C},|s|<\rho(S)\}$. For a matrix $\bld{A}\in \mathbb{R}^{n\times n}$, we use $\sigma^{\rho}(\bld{A})$ and $\sigma^{-}(\bld{A})$ as the shorthand notation for $\sigma^{\rho}(\textup{spec}(\bld{A}))$ and $\sigma^{-}(\textup{spec}(\bld{A}))$, respectively.

A matrix $\bld{A}\in \rn{n}$ is reducible if there exists a permutation matrix \cite{Matrix_Analysis} $\bld{S}\in \rn{n\times n}$ such that $\bld{\hat{A}}=\trn{\bld{S}}\bld{A}\bld{S}=\begin{bmatrix}
\bld{A}_{11} & 0\\
\bld{A}_{21} & \bld{A}_{22}
\end{bmatrix}$. An irreducible matrix is the one that is not reducible. A positive real scalar $p\in \mathbb{R}_{s+}$ is always irreducible.

An irreducible matrix $\bld{A}\in \rn{n\times n}$ is of degree of cyclicity $h$, with $1\leq h\leq n$, if $\sigma^{\rho}(\bld{A})$ is of multiplicity of one with $\sigma^{\rho}(\bld{A})=\{\rho(\bld{A})\textup{exp}(i2\pi k/h),k=0,\dots,h-1\}$  \cite[Th.~2.20]{Plemmons_book}. Moreover, if $\bld{A}\in \rn{n\times n}$ is irreducible with degree of cyclicity $h$, then $\textup{spec}(\bld{A})$ is invariant with respect to polar rotations of $2k\pi/h$ for any $k\in \mathbb{Z}$. 
\subsection{Linear Positive Systems}
\begin{defn}
Define a discrete-time time-invariant linear positive system, with representation
\begin{align}
\bld{x}(t+1)= &\ A\bld{x}(t)+B\bld{u}(t),\ \bld{x(0)}=\bld{x}_0,\label{eq:system_state}\\
\bld{y}(t)= &\ C\bld{x}(t),\label{eq:system_output}
\end{align}
if for any $\bld{x}_{0}\in \rn{n}$ and for any input function $\bld{u}:T\rightarrow\rn{m}$ it holds that the solution of the difference equation \Ref{eq:system_state} is such that $\bld{x}(t)\in \rn{n}$ and $\bld{y}(t)\in \rn{p}$ both for all $t\in T$. Call $A\in \mathbb{R}^{n\times n}$ the system matrix, $B\in \mathbb{R}^{n\times m}$ the input matrix, and $C\in \mathbb{R}^{p\times n}$ the output matrix.  
\end{defn}
It is well known that the solution of the difference equation \Ref{eq:system_state} exists and is provided by the formula,
\begin{equation}\label{eq:system_solution}
\bld{x}(t)=A^{t}\bld{x}_{0}+\sum_{i=1}^{t}A^{i-1}B\bld{u}(t-i).
\end{equation}
Denote this relation by the expression $(0,\bld{x}_0)\overset{\bld{u}(0:t-1)}{\mapsto}(t,\bld{x}(t))$.
\subsection{Controllable Subsets}
\begin{defn}
Consider a discrete-time time-invariant linear positive system with representation 
\begin{align*}
\bld{x}(t+1)= &\ A\bld{x}(t)+B\bld{u}(t),\ \bld{x}(0)=\bld{x}_0,\\
\bld{y}(t)= &\ C\bld{x}(t).
\end{align*}
Define the following subsets of the state space: the $k$-step controllable subset, the finite controllable subset, and the infinite controllable subset, respectively as the sets,
\begin{align}
\cset{k}{\bld{A},\bld{B};\bld{x}_0}&= \{\bld{x}\in \rn{n}|\exists\bld{u}(0:k-1),(0,\bld{x}_0)\mapsto^{\bld{u}(0:k-1)}(k,\bld{x})\},\\ \nonumber
 & k\in \mathbb{Z}^{+},\\ 
\cset{\textup{f}}{\bld{A},\bld{B};\bld{x}_0}&= \cup_{k=0}^{\infty}\cset{\textup{f}}{\bld{A},\bld{B};\bld{x}_0}\\
\cset{\infty}{\bld{A},\bld{B};\bld{x}_0}&= \overline{\cset{\textup{f}}{A,B;\bld{x}_0}}, \forall \bld{x}_0\in \rn{n},
\end{align}
where we have used the notation $\overline{S}$ to denote the closure of the set $S$ with respect to the Euclidean topology. If the initial state equals zero, $\bld{x}_0=0$, then that state is omitted in the notation as in $\cset{k}{\bld{A},\bld{B}}$.
\end{defn}
\section{Characterization of the Controllable Subsets}\label{sect:characterization_Controllable_Subsets}
\begin{prop}\label{prop:controllable_subsets}
Consider a discrete-time linear positive system with the system representation \Ref{eq:system_state} with $\bld{x}_0=0$. The
$k$-step controllable subset, the finite controllable subset, and the infinite controllable subset equal the expressions
\begin{align}
\cset{k}{\bld{A},\bld{B}}&=\textrm{cone}(\textrm{conmat}_{k}(\bld{A},\bld{B})),\label{def:cset_k}\\
\cset{\textup{f}}{\bld{A},\bld{B}}&=\textrm{cone}(\bld{B}\ \bld{AB}\ \bld{A^2B} \ \dots),\label{def:cset_f}\\
\cset{\infty}{\bld{A},\bld{B}}&=\overline{\cset{\textup{f}}{\bld{A},\bld{B}}}, \ \text{where}\label{def:cset_infty}\\
\textrm{conmat}_{k}(\bld{A},\bld{B})&=[\bld{B} \ \bld{AB} \ \bld{A^2B} \ \dots \ \bld{A^{k-1}B}] \label{def:cset_conmat}
\end{align}
\end{prop}
\begin{pf}
Using \Ref{eq:system_solution} with $\bld{x}_0=0$ and with any $\bld{u}:T\rightarrow\rn{m}$, it follows that 
\begin{equation*}
\bld{x}(k)=[\bld{B}\ \bld{AB}\ \dots\ \bld{A}^{k-1}\bld{B}]\trn{[\trn{\bld{u}(k-1)}\ \trn{\bld{u}(k-2)}\ \dots\ \trn{\bld{u}(0)}]}
\end{equation*}
lies in the cone generated by columns of $[\bld{B}\ \bld{AB}\ \dots\ \bld{A}^{k-1}\bld{B}]$ or, equivalently $\cset{k}{\bld{A},\bld{B}}= \textrm{cone}(\textrm{conmat}_k(\bld{A},\bld{B}))$ for any $\bld{u}:T\rightarrow \rn{m}$. The characterization of $\cset{\textup{f}}{\bld{A},\bld{B}}$ and $\cset{\infty}{\bld{A},\bld{B}}$ is then derived in a similar manner.
\end{pf}
\subsection{Polyhedrality of Controllable Subsets}
In this section, given an irreducible matrix $\bld{A}\in \rn{n\times n}$ with degree of cyclicity $1\leq h\leq n$ and $\bld{b}\in \rn{n}$, we first investigate the polyhedrality of $\cset{\infty}{\bld{A},\bld{b}}$, and characterize the necessary and sufficient conditions in terms of $\textup{spec}(\bld{A})$. We then prove that polyhedrality of $\cset{\textup{f}}{\bld{A},\bld{b}}$ is a special case of polyhedrality of $\cset{\infty}{\bld{A},\bld{b}}$ with stricter requirements. In the sequel, it is assumed that $\textup{rank}\big{(}\textup{conmat}_{n}(\bld{A},\bld{b})\big{)}=n$. This condition implies that the characteristic polynomial
and the minimal polynomial coincide \footnote{This is due to the fact that $\bld{A}$ is similar to the companion matrix of $p_{\bld{A}}(\lambda)$ and that for the companion matrix it holds from \cite[pp. 146-147]{Matrix_Analysis} that the characteristic polynomial and the minimal polynomial are equal to $p_{\bld{A}}(\lambda)$.}. This is a convenient assumption that may be relaxed in a future paper.
\begin{prop}\label{prop:polyhedral_conset_infty}
Assume that $A\in \rn{n\times n}$ is irreducible with degree of cyclicity $h\in \mathbb{Z}_{+}$. Define
\begin{align*}
C_{\lim}=&\textup{cone}(\bld{A}_{f,0}\bld{b}, \ldots, \bld{A}_{f,h-1}\bld{b})\\
\bld{A}_{f,i}=&\displaystyle\textup{lim}_{k\rightarrow\infty}\frac{\bld{A}^{kh}}{\rho(\bld{A})^{kh}}\bld{A}^i, \ \textup{for}\ i=0,\dots,h-1.
\end{align*}
Then, the infinite controllable subset $\cset{\infty}{\bld{A},\bld{b}}$ is polyhedral if and only if there exists $k^{\ast}\in\mathbb{Z}_{+}$ such that
\begin{equation}\label{eq:polyhedral_conset_infty_alt_1} 
\textup{cone}(\{\cset{k^{\ast}+1}{\bld{A},\bld{b}}, ~ C_{\lim}\})\subseteq\textup{cone}(\{\cset{k^{\ast}}{\bld{A},\bld{b}}, ~ C_{\lim}\}),
\end{equation}
or equivalently 
\begin{equation}\label{eq:polyhedral_conset_infty_alt} 
\bld{A}^{k^{\ast}}\bld{b}\in \textup{cone}(\{\textup{conmat}_{k^{\ast}}(\bld{A},\bld{b}), ~ C_{\lim}\}).
\end{equation}
\end{prop}
In \Ref{eq:polyhedral_conset_infty_alt_1} and \Ref{eq:polyhedral_conset_infty_alt}, the cone generated by a set of vectors is extended to a cone generated by another cone and a set of vectors. 
\begin{rmk}
Note that due to our assumption on $\bld{A}$, $\bld{A}^h=\textup{diag}(\bld{A}_0,\ldots,\bld{A}_{h-1})$, where $\bld{A_i}\in\rn{n_i\times n_{i}}$, $i=0,\ldots,h-1$, is an irreducible matrix of cyclicity $h=1$ with $\rho(\bld{A}_i)=\rho(\bld{A})^h$, and where $\sum_{i=0}^{h-1} n_{i}=n$. Then, due to \cite[Th.~2.4.1]{Plemmons_book} the limit matrices $\lim_{p\rightarrow \infty}\big{(}\bld{A}_i/\rho(\bld{A}_i)\big{)}^p$, $i=0,\ldots,h-1$ exist. Therefore, the matrices $\bld{A}_{f,i}$ for $i=0,\ldots,h-1$ exist and, hence, the cone $C_{\lim}$ exists.
\end{rmk}
\begin{pf}
Sufficiency: We will show that $$C=\textup{cone}\big{(}\textup{conmat}_{k^{\ast}}(\bld{A},\bld{b})\ \bld{A}_{f,0}\bld{b}\ \dots\ \bld{A}_{f,h-1}\bld{b}\big{)}$$ is $\bld{A}$-invariant. Let $\bld{x}=\sum_{i=0}^{k^{\ast}-1}c_{i}\bld{A}^i\bld{b}+\sum_{i=0}^{h-1}c_{f,i}\bld{A}_{f,i}\bld{b}$ for arbitrary nonnegative coefficients $\bld{c}\in \rn{k^{\ast}}$ and $\bld{c}_{f}\in \rn{h}$. We then have 
\begin{equation}\label{eq:polyhedral_conset_infty_part_1}
\bld{Ax}=\sum_{i=0}^{k^{\ast}-1}c_{i}\bld{A}^{i+1}\bld{b}+\sum_{i=0}^{h-1}c_{f,i}\bld{A}\bld{A}_{f,i}\bld{b}.
\end{equation}
 Using \Ref{eq:polyhedral_conset_infty_alt}, and noting that
\begin{align}
 \bld{A}\bld{A}_{f,i}&=\bld{A}_{f,i+1}, \quad i=0,\dots,h-2\label{eq:A_inf_recursion}\\
\nonumber \bld{A}\bld{A}_{f,h-1}&=\rho(\bld{A})^h\bld{A}_{f,0}, 
\end{align}
\Ref{eq:polyhedral_conset_infty_part_1} can be expressed as $\bld{A}\bld{x}=\sum_{i=0}^{k^{\ast}-1}c^{\prime}_{i}\bld{A}^i\bld{b}+\sum_{i=0}^{h-1}c_{f,i}^{\prime}\bld{A}_{f,i}\bld{b}$ for some $\bld{c}^{\prime}\in \rn{k^{\ast}}$ and some $\bld{c}_{f,i}^{\prime}\in \rn{h}$. This proves $\bld{Ax}\in C$ for any $\bld{x}\in C$. Hence, the system trajectory \Ref{eq:system_solution} remains in $C$ and $\cset{\infty}{\bld{A},\bld{b}}=C$ is polyhedral.

Necessity: Let $\bld{x}_{\infty}=\displaystyle\textup{lim}_{k\rightarrow\infty}\frac{\bld{A}^{k}\bld{b}}{\rho(\bld{A})^{k}}$. Note that even though $\bld{x}_{\infty}$ does not exist in general, its behavior is characterized by the set of $h$ vectors\\ $\bld{A}_{f,0}\bld{b},\dots,\bld{A}_{f,h-1}\bld{b}$ \cite{Coxson1987} (See proof of Lemma~\ref{lem:vf_existence}). Precisely speaking, due to Lemma~\ref{lem:vf_existence}, $\bld{x}_{\infty}\in \textup{cone}(\bld{A}_{f,0}\bld{b},\dots,\bld{A}_{f,h-1}\bld{b})$. By the definition of $\cset{\infty}{\bld{A},\bld{b}}$ as the closure of $\cset{\textup{f}}{\bld{A},\bld{b}}$, and by the above explanation of the vectors $\bld{x}_{\infty}$, the extremal rays of the polyhedral $\cset{\infty}{\bld{A},\bld{b}}$ belong to the sequence $\{ \bld{A}^k \bld{b} \in \mathbb{R}_+, ~ k \in \mathbb{N} \}$
or are extremal rays of the cone, $\textup{cone} (\bld{A}_{f,0}\bld{b}, \ldots, \bld{A}_{f,h-1}\bld{b})$. Again, by the assumption that $\cset{\infty}{\bld{A},\bld{b}}$ is polyhedral, there exists a finite $k^* \in \mathbb{Z}_+$ such that $\bld{A}^{k^*} \bld{b} \in \textup{cone} (\bld{b}, \ldots, \bld{A}^{k^*-1}\bld{b}, \bld{A}_{f,0}\bld{b}, \ldots, \bld{A}_{f,h-1}\bld{b})$.
\end{pf}
It is clear that if \Ref{eq:polyhedral_conset_infty_alt} is established for an integer $k^{\ast}\in \mathbb{Z}_{+}$, it will hold for any $k\geq k^{\ast}$. The smallest integer $k^{\ast}\in \mathbb{Z}_{+}$ satisfying \Ref{eq:polyhedral_conset_infty_alt} is called the \emph{vertex number}, $k_{\textup{vert}}^{\infty}$, of the controllable subset $\cset{\infty}{\bld{A},\bld{b}}$. Following the steps of the proof of Proposition~\ref{prop:polyhedral_conset_infty}, we can put forward the following corollary.
\begin{cor}\label{cor:polyhedral_conset_infty}
The following statements are equivalent:
\begin{enumerate}[(a)]
\item $\cset{\infty}{\bld{A},\bld{b}}$ is polyhedral.
\item There exists an integer $k_{\textup{vert}}^{\infty} \in \mathbb{Z}_+$ such that\\ $\textup{cone}(\bld{b}\ \bld{Ab}\ \dots\ \bld{A}^{k-1}\bld{b}\ \bld{A}_{f,0}\bld{b}\ \dots\ \bld{A}_{f,h-1}\bld{b})$ is $\bld{A}$-invariant for $k\geq k_{\textup{vert}}^{\infty}$.
\item There exists an integer $k_{\textup{vert}}^\infty \in \mathbb{Z}_+$ such that
for the matrix equation,
\begin{eqnarray*}
      \bld{A}\bld{M} 
  & = & \bld{M} \bld{X}, \\
  &   & \exists ~ \mbox{a solution} ~
        \bld{X} \in \mathbb{R}_+^{(k+h) \times (k+h)}, ~ \mbox{with} ~
        k \geq k_{\textup{vert}}^\infty,  ~ \mbox{where,} \\
      \bld{M}
  & = & \left[
        \begin{array}{lllllll}
          \bld{b} & \bld{Ab} & \ldots & \bld{A}^{k-1}\bld{b} & \bld{A}_{f,0}\bld{b} & \ldots & \bld{A}_{f,h-1}\bld{b}
        \end{array}
        \right].
\end{eqnarray*}
\end{enumerate}
\end{cor}
\begin{defn}\label{def:nonnegative_recursion}
A square positive matrix $\bld{A}\in \rn{n\times n}$ is said to have a \emph{nonnegative recursion} if it is satisfied that
\begin{align}
&\exists n_{m}\in \mathbb{N},\ \exists \{c_{0},\dots,c_{n_{m}-1}\}\in \rn{n_{m}}\ \textup{such that}\label{eq:nonnegative_recursion}\\
&\nonumber \bld{A}^{n_{m}}=\sum_{i=0}^{n_{m}-1}c_{i}\bld{A}^i,  
\end{align}
or equivalently
\begin{equation}
g(\lambda)=\lambda^{n_{m}}-\sum_{i=0}^{n_{m}-1}c_{i}\lambda^{i}=0,\ \forall{\lambda}\in \textup{spec}(\bld{A}).
\end{equation}
In terms of the characteristic polynomial, $p_{\bld{A}}(\lambda)$, clearly this implies that 
\begin{equation}
g(\lambda)=p_{\bld{A}}(\lambda)Q(\lambda),
\end{equation}
where $Q(\lambda)$ is a polynomial of degree $n_q\geq 0$. It is immediate that 
\begin{equation}\label{eq:vertex_number_limit}
n_m = n+n_q \geq n.
\end{equation}
\end{defn}
We are now in the position to state a characterization of Proposition~\ref{prop:polyhedral_conset_infty} in terms of spec($\bld{A}$), hence, providing numerically verifiable conditions as to when \Ref{eq:polyhedral_conset_infty_alt} holds. Let
\begin{align*}
\bld{\hat{A}}=\bld{S}^{-1}\bld{A}\bld{S}=
\left[\begin{matrix}
\bld{A}_1 & 0\\
0         & \bld{A}_2
\end{matrix}\right],  
\end{align*}
where $\bld{S}\in \mathbb{R}^{n\times n}$ is non-singular, and where $\bld{A}_1\in \mathbb{R}^{h\times h}$ with $\textup{spec}(\bld{A}_1)=\sigma^{\rho}(\bld{A})$, $\bld{A}_2\in \mathbb{R}^{(n-h)\times (n-h)}$ with $\textup{spec}(\bld{A}_2)=\sigma^{-}(\bld{A})$. Note that such a decomposition is possible due to the Perron-Frobenius theorem \cite[Th.~2.1.4, 2.2.20]{Plemmons_book}. For the pair $(\bld{A},\bld{b})$ of Proposition~\ref{prop:polyhedral_conset_infty} we then have the following theorem.
\begin{thm}\label{thm:polyherdal_conset_infty_characterization}
The following statements are equivalent:
\begin{enumerate}[(a)]
\item \label{item:a_inf} The infinite controllable subset is polyhedral hence there exists an integer $k^* \in \mathbb{Z}_+$ such that
\begin{equation}
\cset{\infty}{\bld{A},\bld{b}}= \textup{cone}(\textup{conmat}_{k^{\ast}}(\bld{A},\bld{b})\ \bld{A}_{f,0}\bld{b}\ \dots\ \bld{A}_{f,h-1}\bld{b}).
\end{equation}
Denote the lowest integer for which the above equality holds by $k_{vert}^{\infty} \in \mathbb{Z}_+$.
\item  \label{item:b_inf} The matrix $\bld{A}_2$ defined above, has a nonnegative recursion.
\item  \label{item:c_inf} If there is a positive $\lambda_{r}\in \textup{spec}(\bld{A}_{2})$, then
\begin{enumerate}[(c1)]
\item $\lambda_r=\rho(\bld{A}_2)$.
\item For any $\lambda\in \sigma^{\rho}(\bld{A}_2)$, $\lambda=\rho(\bld{A}_2)\textup{exp}\big{(}\phi_{\lambda}2\pi i\big{)}$ , where $\phi_{\lambda} \in \mathbb{Q}$ is a rational number.
\item $\sigma^{\rho}(\bld{A}_2)$ are simple.
\item No $\lambda^{-}\in \sigma^{-}(\bld{A}_2)$ has a polar angle which is an integer multiple \footnote{Note that $\sigma^{\rho}(\bld{A}_2)\subseteq \{\lambda\in \mathbb{C}|\lambda=\rho(\bld{A}_2)\textup{exp}\big{(}2k\pi i/(Mh)\big{)}, k=0,\ldots,Mh-1\}$. See Lemma~\ref{lem:subdominant_eigenvalues} for details.} of $2\pi/Mh$.
\end{enumerate} 
\end{enumerate} 
\end{thm}
The proof is established based on a fundamental result \cite[Th.~5]{Roitman1992} on nonnegative recursion, which is also quoted in \ref{sect:Appendix}.
\begin{pf}
\Ref{item:a_inf}$\Rightarrow$\Ref{item:b_inf}$\Rightarrow$\Ref{item:c_inf}:
Since $\cset{\infty}{\bld{A},\bld{b}}$ is polyhedral, according to Corollary~\ref{cor:polyhedral_conset_infty}, there is a sufficiently large $k\geq n-h$ such that the equation    
\begin{equation}
\bld{A}(\bld{b}\ \bld{Ab}\ \dots\ \bld{A}^{k-1}\bld{b}\ \bld{A}_{f,0}\ \dots\ \bld{A}_{f,h-1})=(\bld{b}\ \bld{Ab}\ \dots\ \bld{A}^{k-1}\bld{b}\ \bld{A}_{f,0}\bld{b}\ \dots\ \bld{A}_{f,h-1}\bld{b})\bld{X}
\end{equation}
has a solution $\bld{X}\geq 0$. It can be easily verified using \Ref{eq:polyhedral_conset_infty_alt}-\Ref{eq:A_inf_recursion} that
\begin{equation}\label{eq:X_1}
\bld{X}=\left[\begin{matrix}
\bld{X}_1      & \bld{0}\\
\bld{X}_3      & \bld{X}_2\\
\end{matrix} \right],\
\bld{X_1}=\left[\begin{matrix}
0      & 0      & \cdots & 0 & \alpha_{0}\\
1      & 0      & \cdots & 0 & \alpha_{1}\\
0      & 1      &        & 0 & \alpha_{2}\\
\vdots &        & \ddots &   & \vdots\\
0      & \cdots & 0      & 1 & \alpha_{k}
\end{matrix}\right],
\end{equation}
\begin{equation}\label{eq:X_2}
\bld{X_2}=\left[\begin{matrix}
0      & 0      & \cdots & 0 & \rho(A)^{h}\\
1      & 0      & \cdots & 0 & 0\\
0      & 1      &        & 0 & 0\\
\vdots &        & \ddots &   & \vdots\\
0      & \cdots & 0      & 1 & 0
\end{matrix}\right],\ 
\bld{X_3}=\left[\begin{matrix}
0      & 0      & \cdots & 0 & \beta_{0}\\
0      & 0      & \cdots & 0 & \beta_{1}\\
0      & 0      &        & 0 & \beta_{2}\\
\vdots &        & \ddots &   & \vdots\\
0      & \cdots & 0      & 0 & \beta_{h-1}
\end{matrix}\right].
\end{equation}
constitutes a solution, where $\bld{X}_1\in\rn{k\times k}$, $\bld{X}_2\in\rn{h\times h}$, and $\bld{X}_3\in\rn{h\times k}$. Let $p_{\bld{X}_1}(\lambda)=\det(\lambda \bld{I}-\bld{X}_1)$ and $p_{\bld{X}_2}(\lambda)=\det(\lambda \bld{I}-\bld{X}_2)$. Since, by assumption, $k\geq n-h$ and $\textup{rank}\big{(}\textup{conmat}_{n}(\bld{A},\bld{b})\big{)}=n$, due to \cite[Lemma~3.10]{Benvenuti2006},  $p_{\bld{A}}(\lambda)$ divides $p_{\bld{X}}(\lambda)=p_{\bld{X}_1}(\lambda)p_{\bld{X}_2}(\lambda)=(\lambda^h-\rho(\bld{A})^h)(\lambda^{k}-\alpha_{k-1}\lambda^{k-1}-\dots-\alpha_{0})$. Since $\bld{A}$ is irreducible with degree of cyclicity $h$, $p_{\bld{A}}(\lambda)$ can be expressed as $p_{\bld{A}} (\lambda)=p_{\bld{A}_1}(\lambda)p_{\bld{A}_2}(\lambda)=(\lambda^h-\rho(\bld{A})^h)p_{\bld{A}_2}(\lambda)$. Therefore, $p_{\bld{A}_2}(\lambda)$ divides $p_{\bld{X}_{2}}(\lambda)$, which, due to statements (A) and (B) of Theorem~\ref{thm:Roitman} in \ref{sect:Appendix}, proves $\bld{A}_2$ has a nonnegative recursion of the form $\bld{A}_2^{k^{\ast}}-\alpha_{k^{\ast}-1}\bld{A}_2^{k^{\ast}-1}-\dots-\alpha_{0}\bld{I}=0$ for some $n-h\leq k^{\ast}\leq k$ and for some $\bld{\alpha}\in \rn{k^{\ast}}$. Assume $\bld{A}_2$ has a positive eigenvalue. Since $\bld{A}_2$ satisfies a nonnegative recursion, the statements then (\ref{item:C1}-\ref{item:C4}) in \Ref{item:C} of Theorem~\ref{thm:Roitman} hold for $p_{\bld{A}_2}(\lambda)$. It is straightforward to check that this implies that (\ref{item:c_inf}1)-(\ref{item:c_inf}4) holds\footnote{\label{footnote:equivalence}Condition $\lambda_r\in \sigma^{\rho}(\bld{A}_2)$ follows from \Ref{item:C1} of Theorem~\ref{thm:Roitman}, and conditions (\ref{item:c_inf}2) and (\ref{item:c_inf}3) are, respectively, direct result of \Ref{item:C2} and \Ref{item:C3}. Finally, (\ref{item:c_inf}4) is implied from \Ref{item:C4} using Lemma~\ref{lem:subdominant_eigenvalues}.}. 

\Ref{item:c_inf}$\Rightarrow$\Ref{item:b_inf}$\Rightarrow$\Ref{item:a_inf}: Assume $\bld{A}_2$ has a positive eigenvalue. We need to prove that statements (\ref{item:c_inf}1)-(\ref{item:c_inf}4) imply a nonnegative recursion for $\bld{A}_2$ of the form $\bld{A}_{2}^{k^{\ast}}-\alpha_{k^{\ast}-1}\bld{A}_{2}^{k^{\ast}-1}-\dots-\alpha_{0}\bld{I}=0$, for $k^{\ast}\geq n-h$ and $\bld{\alpha}\in \rn{k^{\ast}}$, and that, in turn, implies polyhedrality of the infinite controllable subset.\\
First we show that the statements (\ref{item:c_inf}1)-(\ref{item:c_inf}4) imply the statements \Ref{item:C1}-\Ref{item:C4} of Theorem~\ref{thm:Roitman}. The statement $\lambda_r\in \sigma^{\rho}(\bld{A}_2)$ implies \Ref{item:C1} of Theorem~\ref{thm:Roitman}. The requirement of all $\lambda\in \sigma^{\rho}(\bld{A}_2)$ having a rational polar phase implies \Ref{item:C2}. The requirement of all $\lambda\in \sigma^{\rho}(\bld{A}_2)$ being simple implies \Ref{item:C3}, and \Ref{item:C4} is implied from $\sigma^{-}(\bld{A}_2)$ including no eigenvalue with polar phase $2\pi m/Mh$ for any $m\in \mathbb{Z}$.
Next, invoking the equivalence between \Ref{item:C} and \Ref{item:B} of Theorem~\ref{thm:Roitman} for $p_{\bld{A_2}}(\lambda)$, one can observe that there is a polynomial $Q(\lambda)$ of positive degree such that 
\begin{equation}\label{eq:thm_2_proof_sufficiency}
g(\lambda)=p_{\bld{A_2}}(\lambda)Q(\lambda)=\lambda^{k^{\ast}}-\alpha_{k^{\ast}-1}\lambda^{k^{\ast}-1}-\dots-\alpha_{0}=0,
\end{equation}
for $k^{\ast}\geq n-h$ and $\bld{\alpha}\in \rn{k^{\ast}}$. It follows from \Ref{eq:nonnegative_recursion} that $\bld{A}_2$ has a nonnegative recursion, which results in \Ref{item:b_inf}.\\ 
Given \Ref{item:b_inf}, there exists a polynomial $g(\lambda)$ of degree $k^{\ast}\geq n-h$ satisfying \Ref{eq:thm_2_proof_sufficiency}, from which one concludes that $p_{\bld{A}}(\lambda)=p_{\bld{A}_1}(\lambda)p_{\bld{A}_2}(\lambda)$ divides $h(\lambda)=p_{\bld{A}_1}(\lambda)g(\lambda)=(\lambda^h-\rho(\bld{A})^h)(\lambda^{k^{\ast}}-\alpha_{k^{\ast}-1}\lambda^{k^{\ast}-1}-\dots-\alpha_{0})$. Now consider the equation $\bld{A}\bld{M}=\bld{M}\bld{X}$ with $\bld{M}=[\bld{b}\ \bld{Ab}\ \dots\ \bld{A}^{k^{\ast}-1}\bld{b}\ \bld{A}_{f,0}\bld{b}\ \dots\ \bld{A}_{f,h-1}\bld{b}]$, where $\bld{X}\in \mathbb{R}^{(n+k^{\ast})\times (n+k^{\ast})}$ is an unknown matrix. Since $\textup{conmat}_{k^{\ast}}(\bld{A},\bld{b})$ is full rank by assumption and $k^{\ast}\geq n-h$, $\bld{M}$ is as well of full rank. Then, it is known from \cite[Lemma~10]{Benvenuti2006} that $p_{\bld{A}}(\lambda)$ divides $p_{\bld{X}}(\lambda)$. Hence, we can choose $\bld{X}$ such that $p_{\bld{X}}(\lambda)=h(\lambda)$. A possible choice of $\bld{X}$, having substituted $k^{\ast}$ for $k$, is then given by \Ref{eq:X_1}-\Ref{eq:X_2}. It is clear from \Ref{eq:X_1}-\Ref{eq:X_2} that $\bld{X}$ admits a nonnegative solution. Based on Corollary~\ref{cor:polyhedral_conset_infty}, this implies that $\cset{\infty}{\bld{A},\bld{b}}$ is polyhedral.
\end{pf}
\begin{rmk}\label{rmk:polyhedral_cset_infty}
For a polyhedral $\cset{\infty}{\bld{A},\bld{b}}$ the following can be observed:
\begin{enumerate}[(a)]
\item Due to \Ref{eq:vertex_number_limit} and from the second part of the proof of Theorem~\ref{thm:polyherdal_conset_infty_characterization} the vertex number of  $\cset{\infty}{\bld{A},\bld{b}}$, $k_{\textup{vert}}^{\infty}$, is at least $n-h$, which implies  $\cset{\infty}{\bld{A},\bld{b}}$ has at least $n$ generators. It has exactly $n$ generators (i.e., is simplicial) if and only if $p_{\bld{A}_2}(\lambda)$ has non-positive coefficients.
\item In the view of Lemma~\ref{lem:vf_existence}, $\cset{\infty}{\bld{A},\bld{b}}$ can be expressed as $\cset{\infty}{\bld{A},\bld{b}}=\textup{cone}(\textup{conmat}_{k_{\textup{vert}}}(\bld{b},
\bld{A})\ \bld{v}_{f,0}\ \dots\ \bld{v}_{f,h-1})$, where $\bld{v}_{f,0},\dots,\bld{v}_{f,h-1}$ are the $h$ distinct nonnegative eigenvectors of $\bld{A}^h$ associated with the eigenvalue $\rho(\bld{A})^h$.
\end{enumerate}
\end{rmk}

\begin{exmp}[polyhedral $\cset{\infty}{\bld{A},\bld{b}}$]\label{ex:1}
Consider the discrete-time linear time-invariant nonnegative system \Ref{eq:system_state} with system matrices
\begin{equation}
\bld{A}=\begin{bmatrix}
0.9727 & 0 & 0.0263\\
0.0388 & 0.1273 & 0.2156\\
0 & 3.4497 & 0
\end{bmatrix},~\bld{b}=\begin{bmatrix}0\\1\\1\end{bmatrix}
\end{equation}
where $\bld{A}$ is primitive, i.e., is irreducible with degree of cyclicity $h=1$. We have $\textup{spec}(\bld{A})=\{1,0.9,-0.8\}$. We can assume $\bld{A}_1=1$, and $\bld{A}_2=\textup{diag}(0.9,-0.8)$. Using Theorem~\ref{thm:polyherdal_conset_infty_characterization}, it is immediate that conditions (\ref{item:c_inf}1) and (\ref{item:c_inf}2) hold as $\lambda=0.9$ is a simple eigenvalue of $\bld{A}_2$, which equals the spectral radius of $\bld{A}_2$. Condition (\ref{item:c_inf}1) hold as well since the polar angle of $\lambda=-0.8$ is not a integer multiple of the polar angle of $\lambda=0.9$. Hence, it can be concluded that the infinite controllable subset $\cset{\infty}{\bld{A},\bld{b}}$ is polyhedral. We can also conclude that $\bld{A}_2$ has a nonnegative recursion, which is readily verified as $p_{\bld{A}_2}(\lambda)=\lambda^2-0.1\lambda-0.72$. Fig.~\ref{fig:E1} illustrates the growth of $\cset{k}{\bld{A},\bld{b}}$. It can be observed that $\cset{\textup{f}}{\bld{A},\bld{b}}$ is not polyhedral since the cone keeps growing for increasing values of $k$. Its closure is, however, polyhedral as shown in Fig.~\ref{fig:E1}d.
\begin{figure}
\hspace*{\fill}
\begin{minipage}[b]{0.45\linewidth}
\centering
\includegraphics[width=\textwidth]{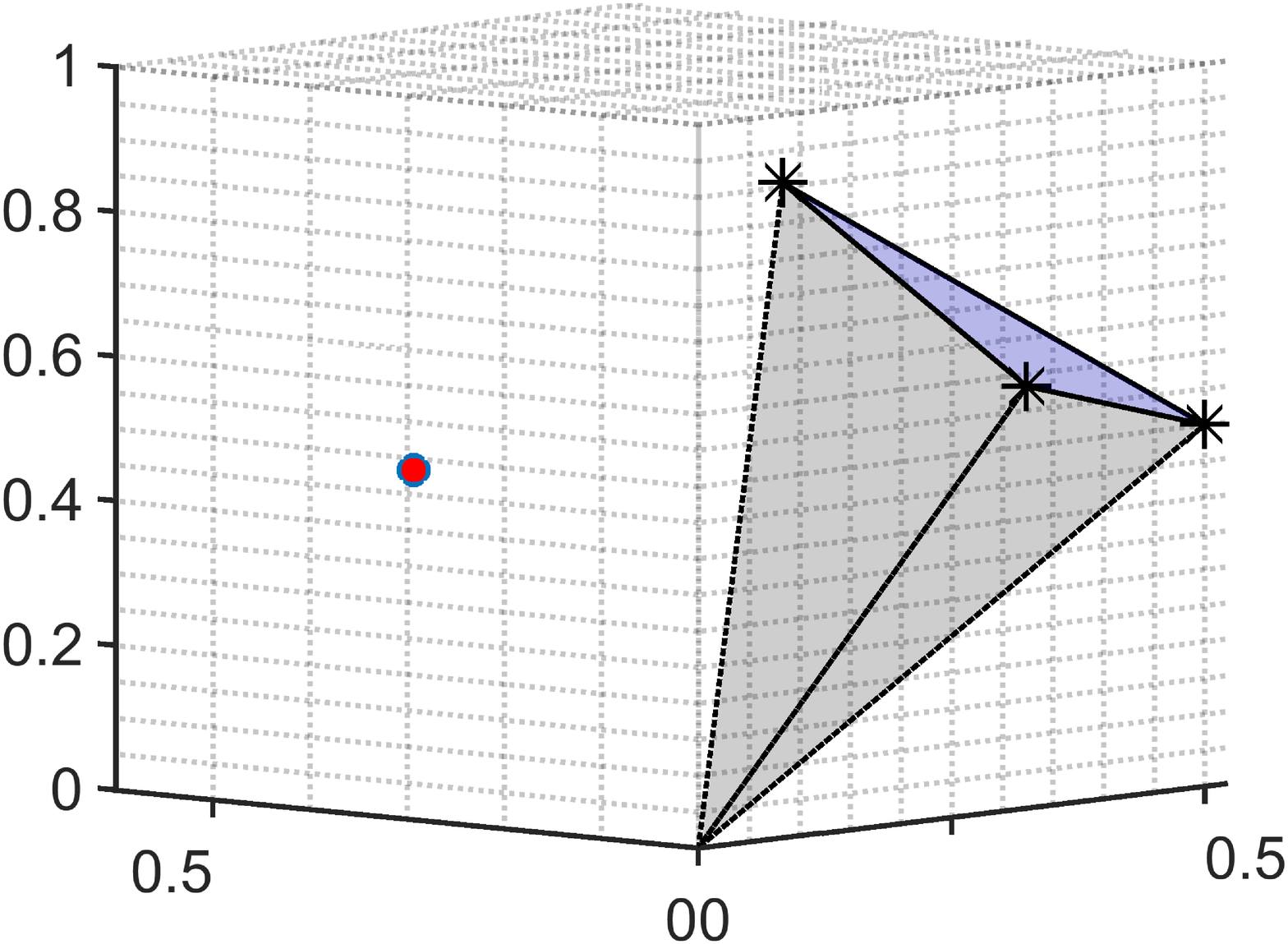}
\captionsetup{labelformat=empty}
\caption*{(a)~$\cset{3}{\bld{A},\bld{b}}$}
\end{minipage}
\hspace*{\fill}
\begin{minipage}[b]{0.45\linewidth}
\centering
\includegraphics[width=\textwidth]{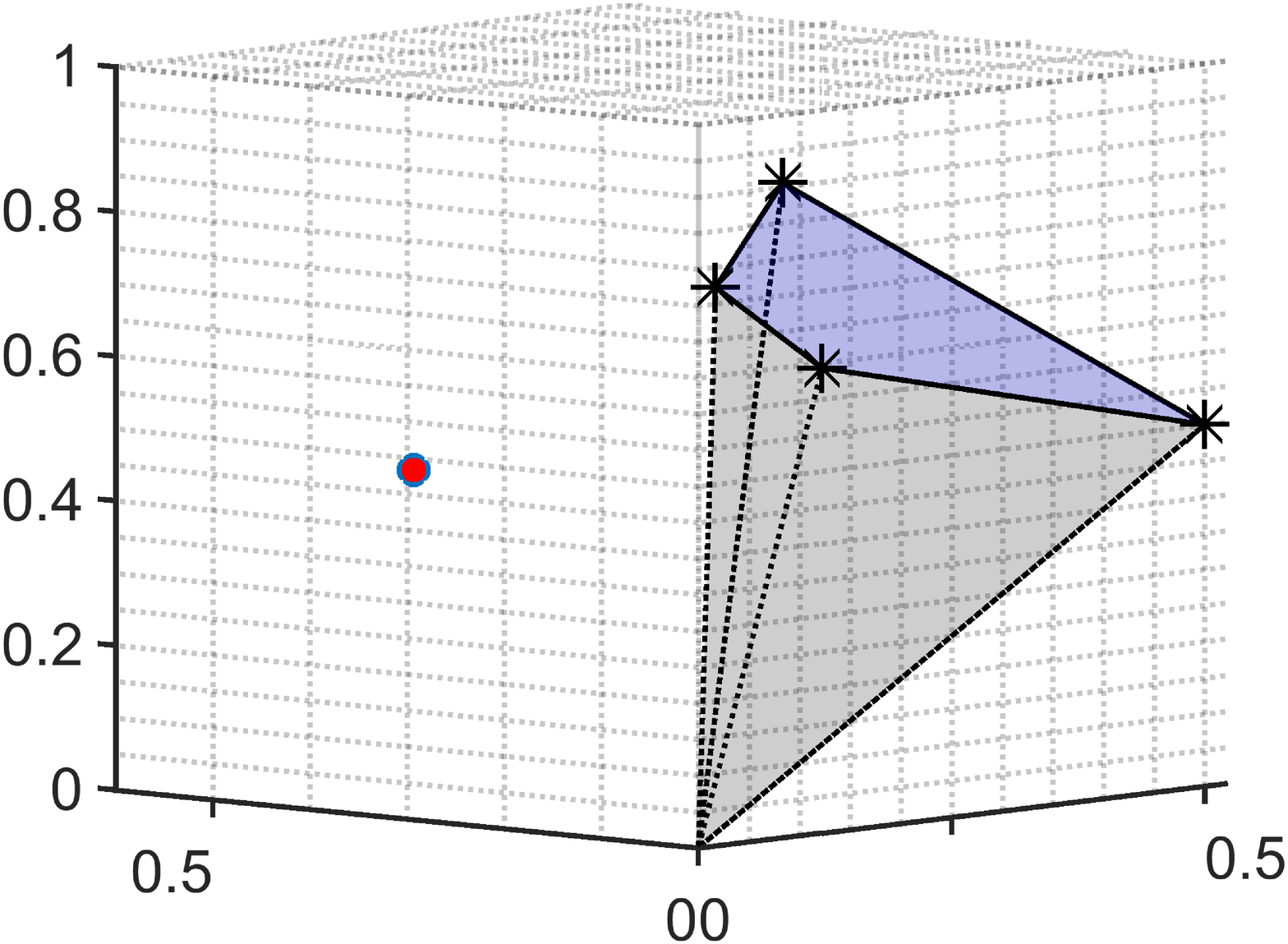}
\captionsetup{labelformat=empty}
\caption*{(b)~$\cset{8}{\bld{A},\bld{b}}$}
\end{minipage}
\hspace*{\fill}
\\
\hspace*{\fill}
\begin{minipage}[b]{0.45\linewidth}
\centering
\includegraphics[width=\textwidth]{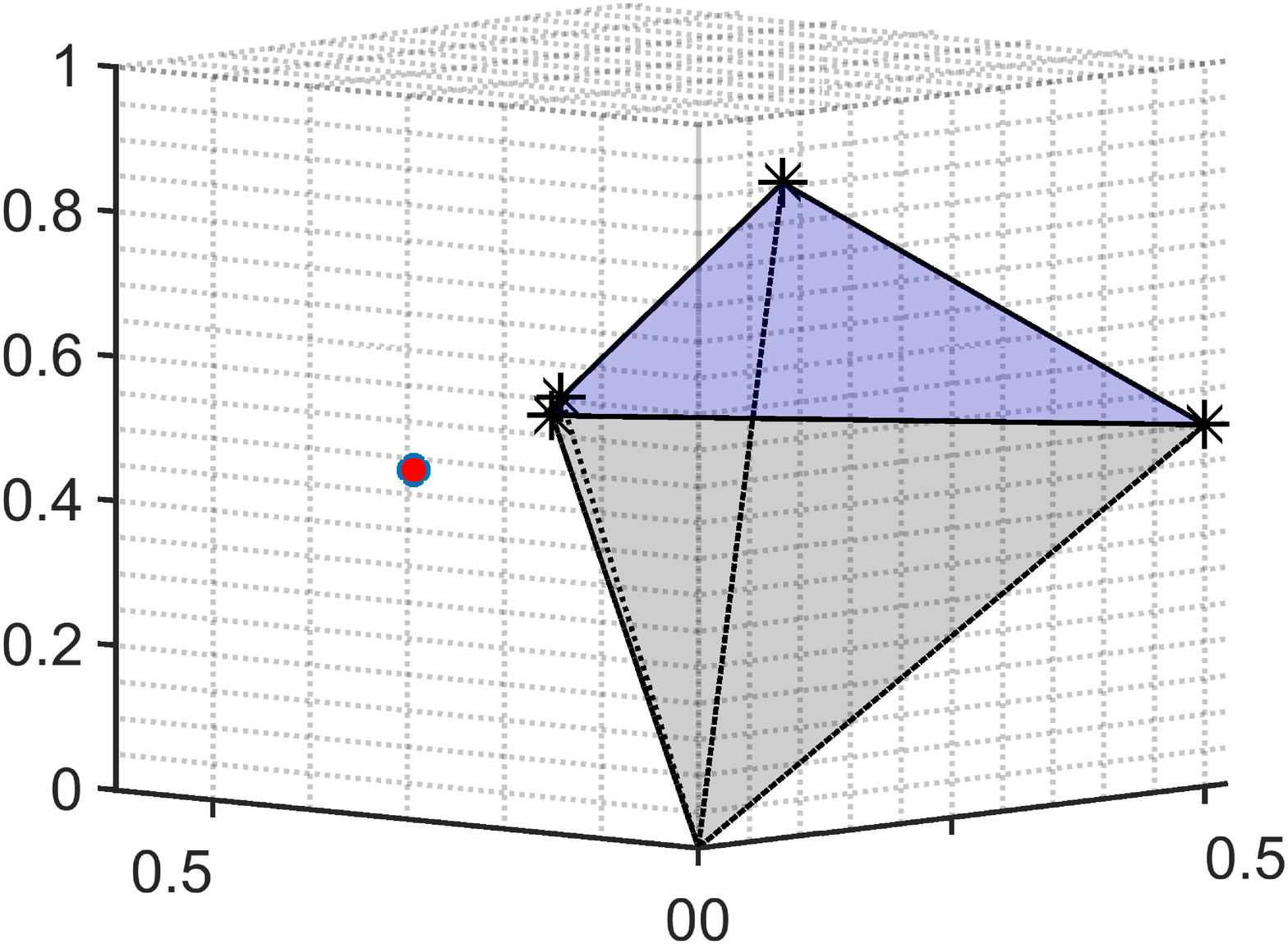}
\captionsetup{labelformat=empty}
\caption*{(c)~$\cset{19}{\bld{A},\bld{b}}$}
\end{minipage}
\hspace*{\fill}
\begin{minipage}[b]{0.45\linewidth}
\centering
\includegraphics[width=\textwidth]{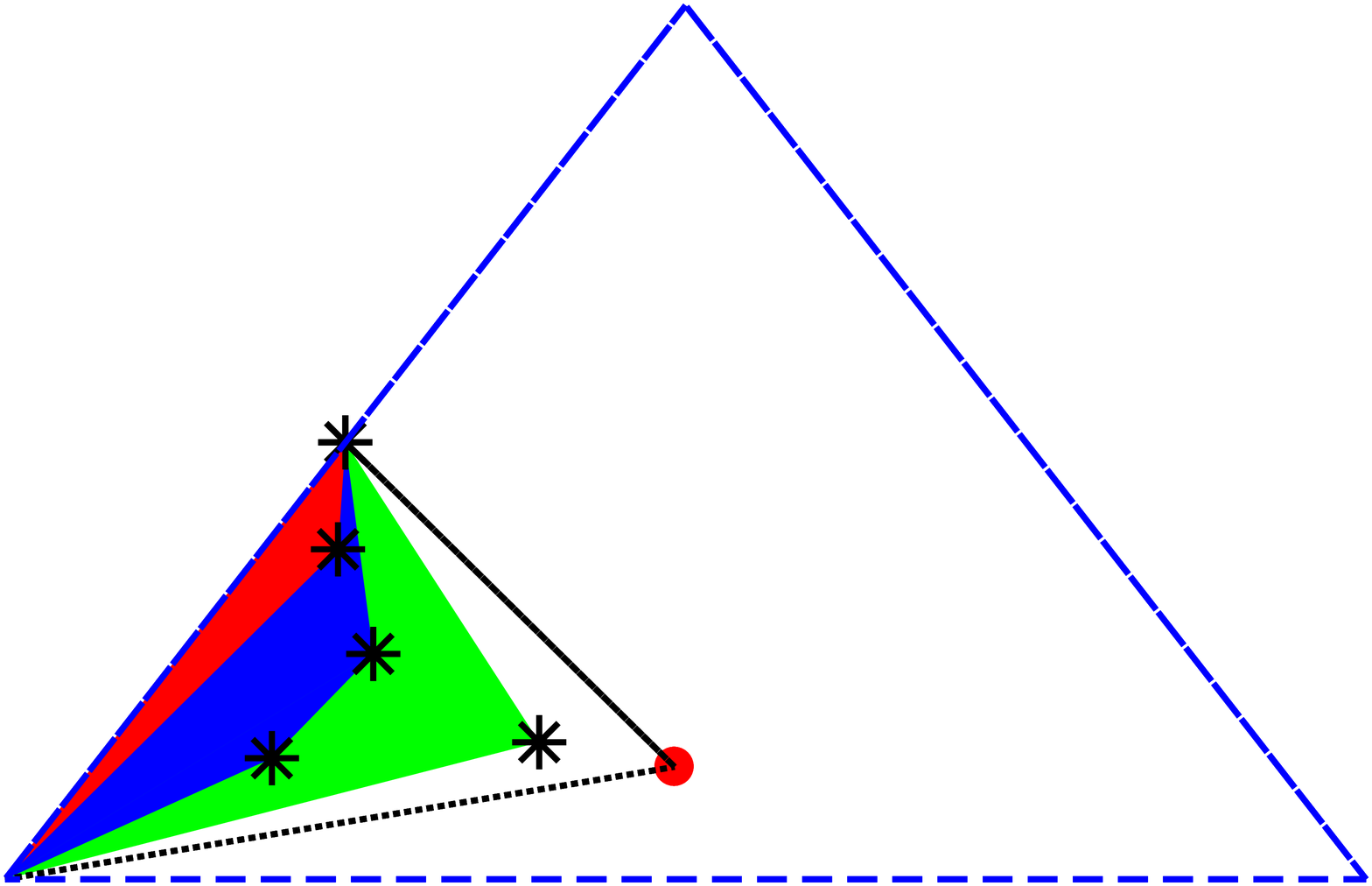}
\captionsetup{labelformat=empty}
\caption*{(d)~$\cset{k}{\bld{A},\bld{b}}$, $k$=3 (red), 8 (blue and red), 19 (green, blue and red) and $\cset{\infty}{\bld{A},\bld{b}}$(the triangle with the red vertex)}
\end{minipage}
\hspace*{\fill}
\caption{\label{fig:E1}a,b,c: the growth of controllable cone $\cset{k}{\bld{A},\bld{b}}$ of example~\ref{ex:1} for different values of $k$, where generators of the cone are marked by asterisks, and the Frobenius eigenvector is marked by a red dot. d: The growth of controllable cone mapped on the 3-dimensional simplex $S=\{\bld{x}\in \mathbb{R}^{3}_{+}|\trn{\mathbbm{1}}\bld{x}=1\}$.}
\end{figure}
\end{exmp}
\begin{exmp}[non-polyhedral $\cset{\infty}{\bld{A},\bld{b}}$]\label{ex:2}
Consider the discrete-time linear time-invariant nonnegative system \Ref{eq:system_state} with system matrices
\begin{equation}
\bld{A}=\begin{bmatrix}
0 & 1 & 0\\
1 & 0 & 0.5\\
0 & 0.4 & 1
\end{bmatrix},~\bld{b}=\begin{bmatrix}0\\1\\0\end{bmatrix},
\end{equation}
where $\bld{A}$ has degree of cyclicity $h=1$. The spectrum of $\bld{A}$ is $\textup{spec}(\bld{A})=\{-1.05,0.7116,1.3383\}$. One can assume $\bld{A}_1=1.3383$ and\\ $\bld{A}_2=\textup{diag}(-1.05,0.7116)$. It is immediate that condition (c1) of Theorem~\ref{thm:polyherdal_conset_infty_characterization} is not satisfied as $0.7116\neq \rho(\bld{A}_2)$. Therefore, based on this theorem, $\cset{\infty}{\bld{A},\bld{b}}$ is not polyhedral. This is illustrated by Fig.~\ref{fig:E2}d, from which it is clear that $\cset{\infty}{\bld{A},\bld{b}}$ is approaching a round cone as introduced in Section~\ref{subsect:Positive_Marices}.
\begin{figure}
\hspace*{\fill}
\begin{minipage}[b]{0.45\linewidth}
\centering
\includegraphics[width=\textwidth]{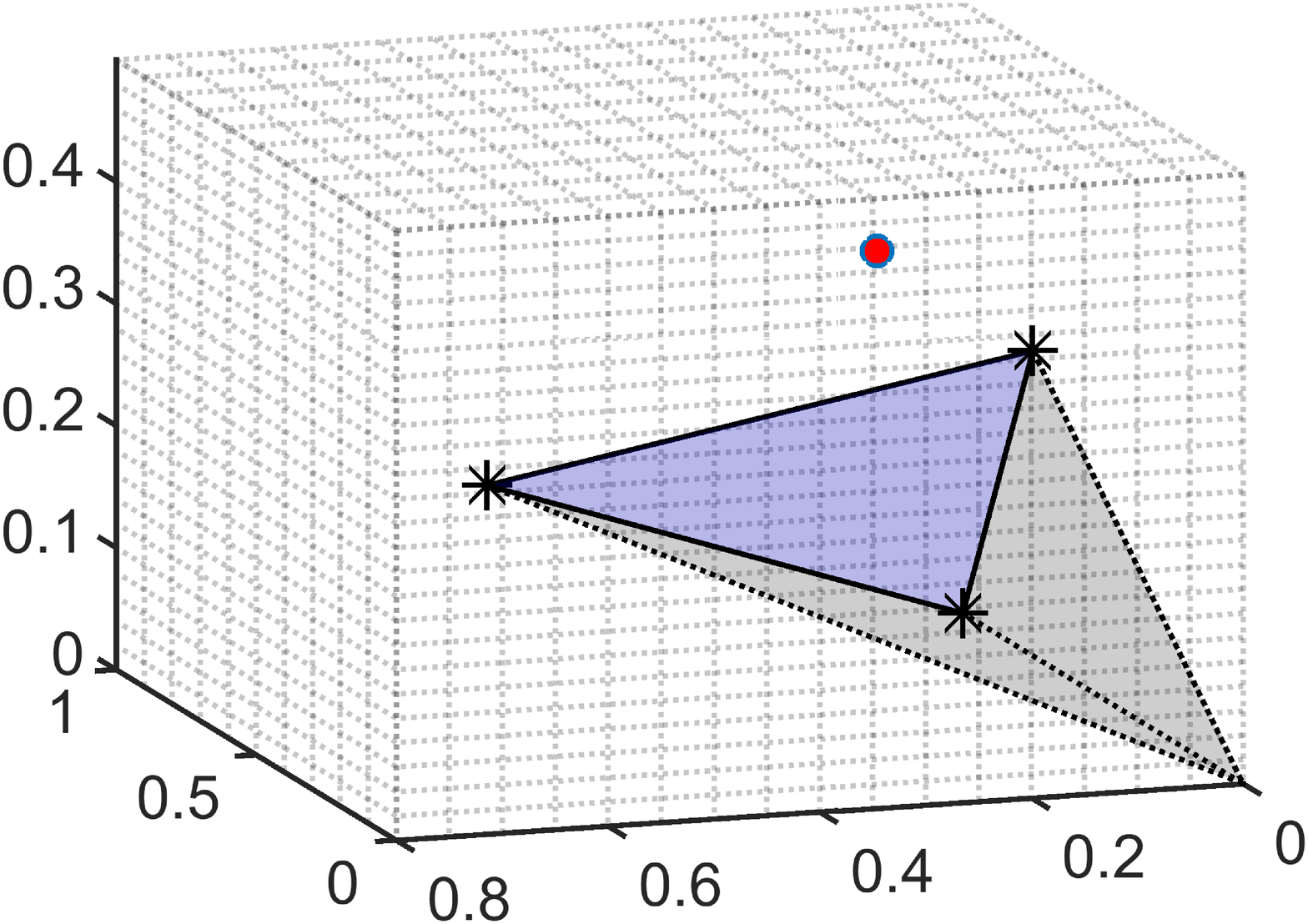}
\captionsetup{labelformat=empty}
\caption*{(a)~$\cset{k}{\bld{A},\bld{b}},~k=3$}
\end{minipage}
\hspace*{\fill}
\begin{minipage}[b]{0.45\linewidth}
\centering
\includegraphics[width=\textwidth]{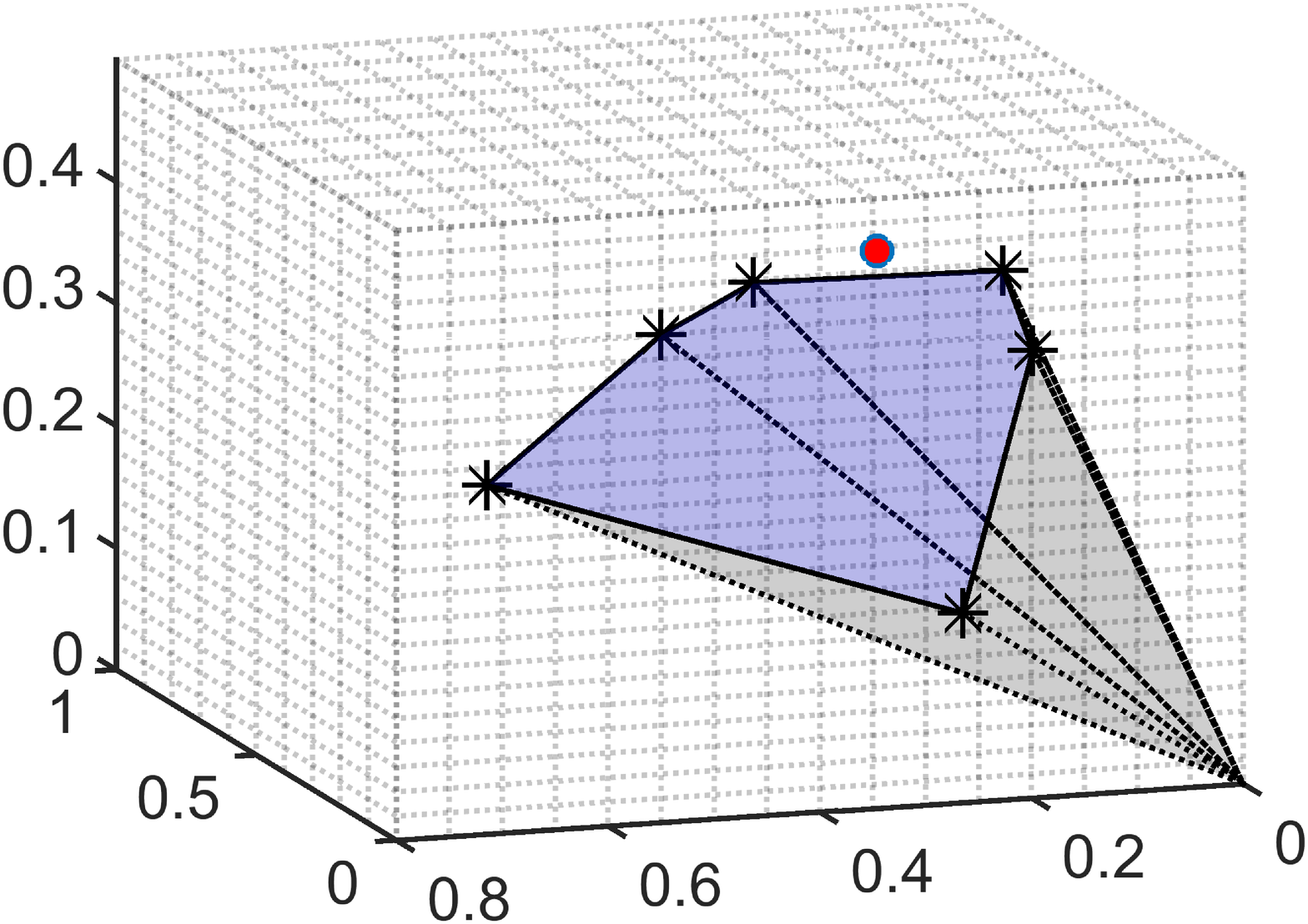}
\captionsetup{labelformat=empty}
\caption*{(b)~$\cset{k}{\bld{A},\bld{b}},~k=6$}
\end{minipage}
\hspace*{\fill}
\\
\hspace*{\fill}
\begin{minipage}[b]{0.45\linewidth}
\centering
\includegraphics[width=\textwidth]{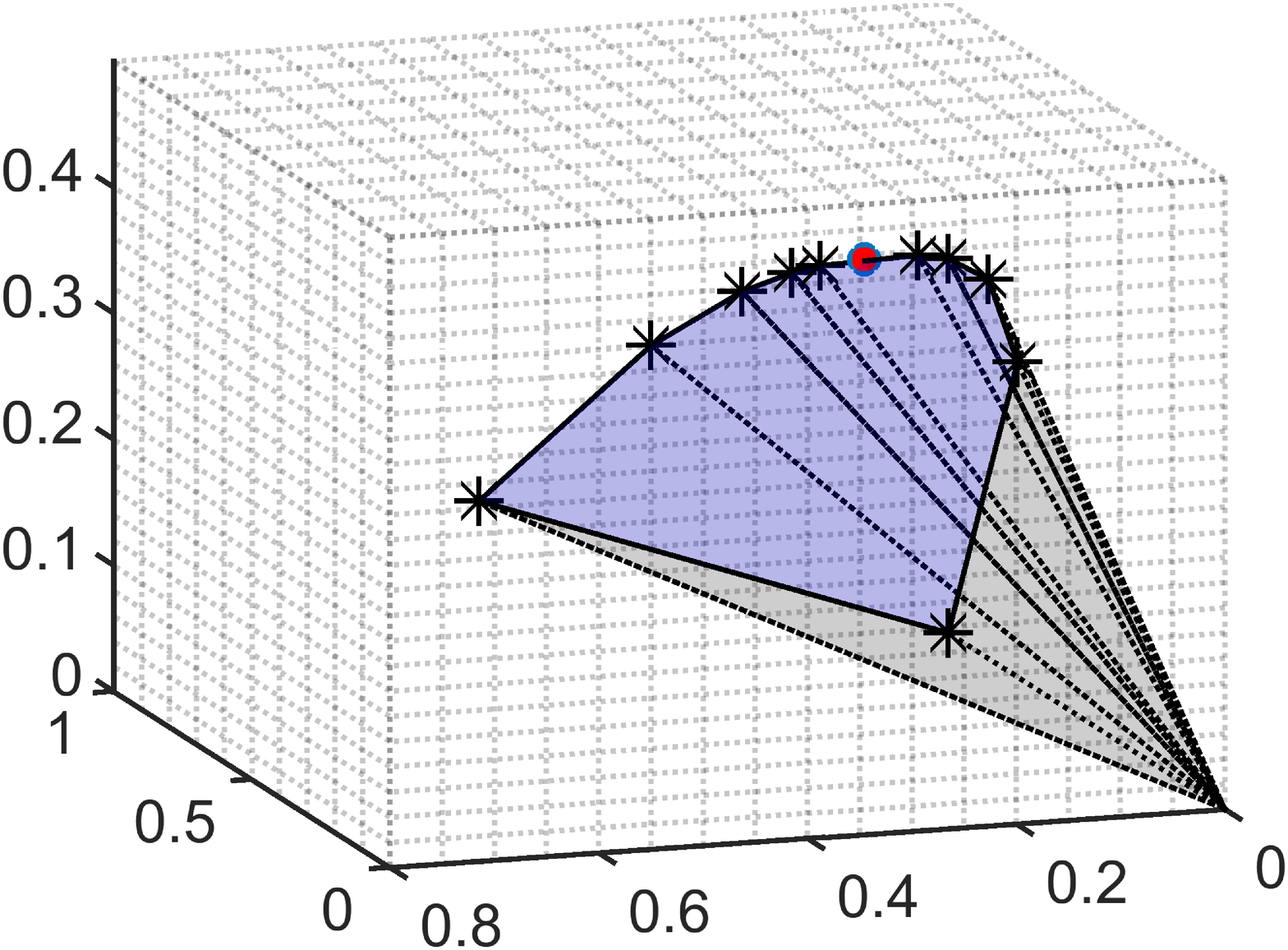}
\captionsetup{labelformat=empty}
\caption*{(c)~$\cset{k}{\bld{A},\bld{b}},~k=10$}
\end{minipage}
\hspace*{\fill}
\begin{minipage}[b]{0.45\linewidth}
\centering
\includegraphics[width=\textwidth]{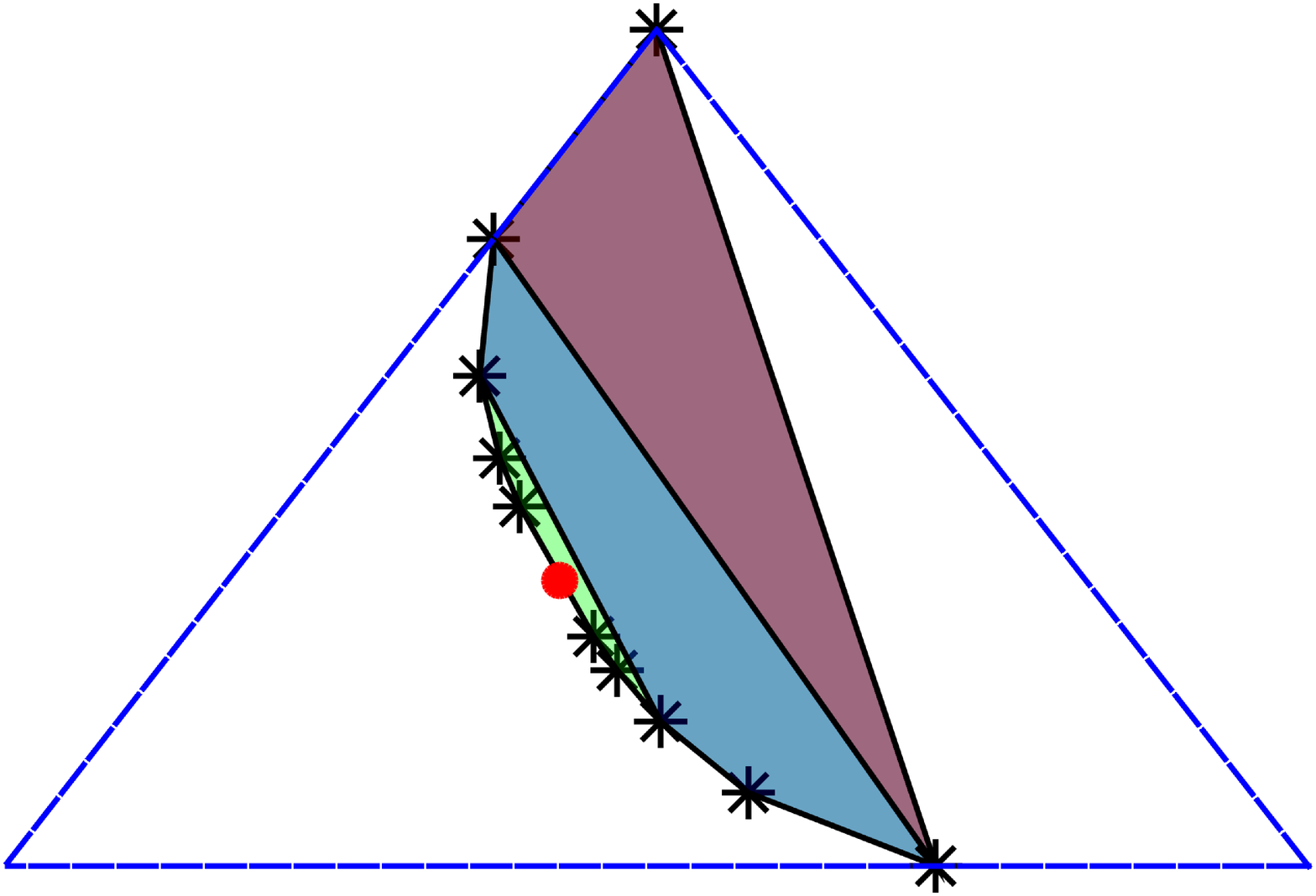}
\captionsetup{labelformat=empty}
\caption*{(d)~$\cset{k}{\bld{A},\bld{b}}$, $k=3$ (red region), $k=6$ (red and blue regions), $k=10$ (red, blue and green regions). $\cset{\infty}{\bld{A},\bld{b}}$ approaches a ``round cone''.}
\end{minipage}
\hspace*{\fill}
\caption{\label{fig:E2}a,b,c: the growth of controllable cone $\cset{k}{\bld{A},\bld{b}}$ of example~\ref{ex:2} for different values of $k$, where generators of the cone are marked by asterisks, and the Frobenius eigenvector is marked by a red dot. d: The growth of controllable cone mapped on the 3-dimensional simplex $S=\{\bld{x}\in \mathbb{R}^{3}_{+}|\trn{\mathbbm{1}}\bld{x}=1\}$.}
\end{figure}
\end{exmp}
Now we will investigate polyhedrality of $\cset{\textup{f}}{\bld{A},\bld{b}}$. Consider the following proposition. We will show that this implies stricter conditions on $\textup{spec}(\bld{A})$ and that a more conservative version of Theorem~\ref{thm:polyherdal_conset_infty_characterization} applies. 
\begin{prop}\label{prop:polyhedral_conset_f}
The finite controllable subset $\cset{\textup{f}}{\bld{A},\bld{b}}$ is polyhedral if and only if there exists a positive integer $k^{\ast}\in \mathbb{Z}_{+}$ such that
\begin{equation}
\cset{k^{\ast}+1}{\bld{A},\bld{b}} \subseteq \cset{k^{\ast}}{\bld{A},\bld{b}}, 
\end{equation}
or equivalently,
\begin{equation}\label{eq:polyhedral_conset_f}
\bld{A}^{k^{\ast}}\bld{b}\in \cset{k^{\ast}}{\bld{A},\bld{b}}.
\end{equation}
\end{prop}
\begin{pf}
Sufficiency: If $\bld{A}^{k^{\ast}}\bld{b}\in \cset{k^{\ast}}{\bld{A},\bld{b}}$ it follows immediately from \Ref{eq:system_solution} that $\bld{x}(t)\in \cset{k^{\ast}}{\bld{A},\bld{b}}$ for any $t\geq k^{\ast}$. Hence, based on \Ref{def:cset_f} in Proposition~\ref{prop:controllable_subsets}, $\cset{\textup{f}}{\bld{A},\bld{b}}=\cset{k^{\ast}}{\bld{A},\bld{b}}$.\\
Necessity: if $\cset{\textup{f}}{\bld{A},\bld{b}}$ is a polyhedral cone, since its generators are of the form $\bld{A}^k\bld{b}$, $k\in \mathbb{Z}_{+}$,  and since a polyhedral cone has a finite number of generators, there must exist a finite $k^{\ast}\in \mathbb{Z}_{+}$ for which $\bld{A}^{k^{\ast}}\bld{b}\in \cset{k^{\ast}}{\bld{A},\bld{b}}$.
\end{pf}
The smallest $k^{\ast}$ for which \Ref{eq:polyhedral_conset_f} holds is referred to as the \emph{vertex number}, $k_{\textup{vert}}$, of  $\cset{\textup{f}}{\bld{A},\bld{b}}$. Note that \Ref{eq:polyhedral_conset_f} also implies that
\begin{equation}\label{eq:cset_f_polyhedrality_restriction}
\textup{cone}(\bld{A}_{f,0}\bld{b}\ \dots\ \bld{A}_{f,h-1}\bld{b})\subset\cset{k_{\textup{vert}}}{\bld{A},\bld{b}},
\end{equation}
which is clearly a restriction on \Ref{eq:polyhedral_conset_infty_alt}. Based on the proof of Proposition~\ref{prop:polyhedral_conset_f}, one can derive the following corollary.
\begin{cor}\label{cor:polyhedral_conset_f}
The following statements regarding polyhedrality of $\cset{\textup{f}}{\bld{A},\bld{b}}$  are equivalent:
\begin{enumerate}[(a)]
\item $\cset{\textup{f}}{\bld{A},\bld{b}}$ is polyhedral.
\item There exists an integer $k_{\textup{vert}} \in \mathbb{Z}_+$ such that $\textup{cone}(\bld{b}\ \bld{Ab}\ \dots\ \bld{A}^{k}\bld{b})$ is $\bld{A}$-invariant for any $k\geq k_{\textup{vert}}$.
\item There exists an integer $k_{\textup{vert}} \in \mathbb{Z}_+$ such that for the matrix equation,
\begin{eqnarray*}
      \bld{A}[\bld{b}\ \bld{Ab}\ \dots\ \bld{A}^{k-1}\bld{b}]&=&[\bld{b}\ \bld{Ab}\ \dots\ \bld{A}^{k-1}\bld{b}]\bld{X},\\
  &   & \exists ~ \mbox{a solution} ~
        \bld{X} \in \mathbb{R}_+^{(k) \times (k)}, ~ \mbox{with} ~
        k \geq k_{\textup{vert}}.
\end{eqnarray*}
\item Based on \Ref{eq:cset_f_polyhedrality_restriction} and Lemma~\ref{lem:vf_existence}, there exists an integer $k_{\textup{vert}} \in \mathbb{Z}_+$ such that $\textup{cone}(\bld{v}_{f,0}\ \dots\ \bld{v}_{f,h-1})\subset\cset{k}{\bld{A},\bld{b}}$ for any $k\geq k_{\textup{vert}}$. 
\end{enumerate}
\end{cor}
Now, a decomposition of $\bld{A}$ is introduced that will be used for stating the next theorem. Given $\bld{A}\in \rn{n\times n}$, consider $\bld{A}_1\in \mathbb{R}$ and $\bld{A}_2\in \mathbb{R}^{(n-1)\times (n-1)}$, where $\textup{spec}(\bld{A}_1)=\rho(\bld{A})$ and $\textup{spec}(\bld{A}_2)=\textup{spec}(\bld{A})\setminus\{\rho(\bld{A})\}$. The decomposition of $\bld{A}$ into $\bld{A}_1$ and $\bld{A}_2$ is then given by $\bld{A}=\bld{S}\textup{diag}(\bld{A}_1,\bld{A_2})\bld{S}^{-1}$, where $\bld{S}\in \mathbb{R}^{n\times n}$ is non-singular. Note that such a decomposition is possible due to the Perron-Frobenius theorem \cite[Th.~2.1.4, 2.2.20]{Plemmons_book}. With such decomposition of $\bld{A}$ at hand, the following theorem provides necessary and sufficient conditions on $\textup{spec}(\bld{A})$ for polyhedrality of $\cset{\textup{f}}{\bld{A},\bld{b}}$. These conditions turn out to be a conservative version of those of Theorem~\ref{thm:polyherdal_conset_infty_characterization}. 
\begin{thm}\label{thm:polyherdal_conset_f_characterization}
The following statements are equivalent:
\begin{enumerate}[(a)]
\item \label{item:a}The finite controllable subset is polyhedral and hence there exists an integer $k^{\ast} \in \mathbb{Z}_+$, $k^{\ast}\geq k_{\textup{vert}}$ such that $\cset{\textup{f}}{\bld{A},\bld{b}}=\cset{k^{\ast}}{\bld{A},\bld{b}}$.
\item \label{item:b}$\bld{A}$ has a nonnegative recursion.
\item \label{item:c}$\bld{A}_2$ does not have any positive eigenvalue.  
\end{enumerate}
\end{thm}
\begin{pf}
\Ref{item:a} $\Rightarrow$ \Ref{item:b} $\Rightarrow$ \Ref{item:c}: Based on Corollary~\ref{cor:polyhedral_conset_f} with $k\geq n$ we obtain 
$$\bld{A}\big{(}\textup{conmat}_{k}(\bld{A},\bld{b})\big{)}=\big{(}\textup{conmat}_{k}(\bld{A},\bld{b})\big{)}\bld{X},$$
where $\bld{X}\in \rn{k\times k}$ is given by 
$$\bld{X}=\left[\begin{matrix}
0      & 0      & \cdots & 0 & \alpha_{0}\\
1      & 0      & \cdots & 0 & \alpha_{1}\\
0      & 1      &        & 0 & \alpha_{2}\\
\vdots &        & \ddots &   & \vdots\\
0      & \cdots & 0      & 1 & \alpha_{k-1}
\end{matrix}\right].$$ Since, by assumption, $\textup{conmat}_{n}(\bld{A},\bld{b})$ is full rank and $k\geq n$, there exists \cite[Lemma~3.10]{Benvenuti2006} a polynomial $Q(\lambda)$ of nonnegative degree such that $p_{\bld{A}}(\lambda)Q(\lambda)=p_{\bld{X}}(\lambda)=\poly{\lambda}{k}$, which, in the view of Definition~\ref{def:nonnegative_recursion}, proves that $\bld{A}$ has a nonnegative recursion. Noting that \Ref{item:b} is equivalent to condition \Ref{item:B} of Theorem~\ref{thm:Roitman} (\cite[Th.~5]{Roitman1992}), all conditions \Ref{item:C1}-\Ref{item:C4} are then fulfilled. In particular, \Ref{item:C4} holds as conditions \Ref{item:C1}-\Ref{item:C3} are already satisfied for a nonnegative irreducible matrix due to the Perron-Frobenius theorem \cite[Th.~2.1.4, 2.2.20]{Plemmons_book}. Condition \Ref{item:C4} requires that no eigenvalue $\lambda^{-}\in \sigma^{-}(\bld{A})$ has a polar angle of $2\pi k/h$ for $k=0,\dots,h-1$. Since $\textup{spec}(\bld{A})$ is invariant under a polar rotation of $2\pi m/h$ for any $m\in \mathbb{Z}$, no $\lambda^{-}\in \sigma^{-}(\bld{A})$ is then positive. Noting that for an irreducible matrix,  $\big{(}\sigma^{\rho}(\bld{A})\setminus\{\rho(\bld{A})\}\big{)}\cap \mathbb{R}_+=\emptyset$ and that $\textup{spec}(\bld{A}_2)=\sigma^{-}(\bld{A})\cup\sigma^{\rho}(\bld{A})\setminus\{\rho(\bld{A})\}$, one concludes that $\bld{A}_2$ has no positive eigenvalue.

\Ref{item:c} $\Rightarrow$ \Ref{item:b} $\Rightarrow$ \Ref{item:a}: Given \Ref{item:c}, we have $\textup{spec}(\bld{A}_2)\cap \mathbb{R}_{+}=\emptyset$. For an irreducible matrix it holds that $\big{(}\sigma^{\rho}(\bld{A})\setminus\{\rho(\bld{A})\}\big{)}\cap \mathbb{R}_+=\emptyset$. Since $\textup{spec}(\bld{A}_2)=\sigma^{-}(\bld{A})\cup (\sigma^{\rho}(\bld{A})\setminus\{\rho(\bld{A})\})$, it follows that $\sigma^{-}(\bld{A})\cap \mathbb{R}_{+}=\emptyset$, from which it can be immediately concluded that $\not\exists \lambda\in \sigma^{-}(\bld{A}), ~\lambda=|\lambda|\textup{exp}(i2\pi m/h)$ for any $m\in \mathbb{Z}$. Hence, we establised that \Ref{item:C4} of Theorem~\ref{thm:Roitman} (\cite[Th.~5]{Roitman1992}) holds for $p_{\bld{A}}(\lambda)$. Moreover, statements \Ref{item:C1}-\Ref{item:C3} as well hold for $p_{\bld{A}}(\lambda)$ as $\bld{A}$ is irreducible. Therefore, due to \Ref{item:B} of Theorem \ref{thm:Roitman}, there exists a polynomial $Q(\lambda)$ of nonnegative degree, such that $p_{\bld{A}}(\lambda)Q(\lambda)=\poly{\lambda}{k^{\ast}}$, where $k^{\ast}\geq n$ and $\alpha_{i}\geq 0$, $i=0,1,\dots,k^{\ast}-1$. This proves that $\bld{A}$ has a nonnegative recursion based on Definition \ref{def:nonnegative_recursion}. Then, \Ref{item:a} immediately follows as $\bld{A}^{k^{\ast}}\bld{b}=\sum_{i=0}^{k^{\ast}-1}\alpha_{i}\bld{A}^{i}\bld{b}$. 
\end{pf}
\begin{rmk}\label{rmk:polyhedral_cset_f}
Note that since $\textup{deg}\big{(}Q(\lambda)\big{)}\geq 0$, $k_{\textup{vert}}$ of $\cset{\textup{f}}{\bld{A},\bld{b}}$ is at least $n$, and it equals $n$ if and only if $p_{\bld{A}}(\lambda)=\poly{\lambda}{n}$ with $\alpha_{i}\geq 0$, $i=0,\dots,n-1$. Hence $\cset{\textup{f}}{\bld{A},\bld{b}}$ is a simplicial cone (i.e., has $n$ generators) if and only if the characteristic polynomial of $\bld{A}$ has non-positive coefficients. One such matrix is a cyclic matrix with cyclicity index $h=n$ as $p_{\bld{A}}(\lambda)=\lambda^n-\rho(\bld{A})^n$. 
\end{rmk}
Comparing Theorem~\ref{thm:polyherdal_conset_infty_characterization} to Theorem~\ref{thm:polyherdal_conset_f_characterization} reveals that the latter is a restricted version of the former. For example, Theorem~\ref{thm:polyherdal_conset_infty_characterization}b requires a part of $\bld{A}$($\bld{A}_2$) to have a nonnegative recursion while Theorem~\ref{thm:polyherdal_conset_f_characterization}b requires $\bld{A}$ to have a nonnegative recursion.
\begin{exmp}[polyhedral $\cset{\textup{f}}{\bld{A},\bld{b}}$]\label{ex:3}
Consider the discrete-time linear \\time-invariant nonnegative system \Ref{eq:system_state} with system matrices
\begin{equation}
\bld{A}=\begin{bmatrix} 
0 & 1.6333 & 1.1049 & 0\\
23.5667 & 6.0944 & 0 & 0\\
0 & 0 & 1.1225 & 1.0672\\
0 & 1.6611 & 0 & 0.7830
\end{bmatrix},\ \bld{b}=\begin{bmatrix}0\\0\\1\\1 \end{bmatrix},
\end{equation}
where $\bld{A}$ is irreducible with degree of cyclicity $h=1$. It can be verified that $\textup{spec}(\bld{A})=\{10,-4,1+1i,1-1i\}$. One can recognize that no eigenvalue of $\bld{A}_2=\textup{diag}(-4,1+i,1-i)$ is positive. Therefore, condition (c3) of Theorem~\ref{thm:polyherdal_conset_f_characterization} holds and it follows that $\bld{A}$ has a nonnegative recursion. In fact, it can be verified that in this case it holds that $\bld{A}^6 = 166.7569\bld{I}_{4}+16.1434\bld{A}+39.7036\bld{A}^4+6.0262\bld{A}^5$, where $\bld{I}_4$ denotes the identity matrix of dimension $4\times 4$. In addition, we can conclude that $\cset{\textup{f}}{\bld{A},\bld{b}}$ is polyhedral with $k_{\textup{vert}}=6$. This is illustrated by Fig.~\ref{fig:E3}, where it is observed that $\cset{k}{\bld{A},\bld{b}}$ stops growing for $k\geq 6$, that is $\cset{k}{\bld{A},\bld{b}}=\cset{6}{\bld{A},\bld{b}}$ for any $k\geq 6$. One can also notice from Fig.~\ref{fig:E3}c that $C_{\lim}\subset \cset{k_{\textup{vert}}}{\bld{A},\bld{b}}$. Note that in this particular example, since $h=1$, we have $C_{\lim}=\textup{cone}(\bld{A}_{f,0}\bld{b})=\{c\bld{v}_f|c \in \mathbb{R}_{+}\}$, where $\bld{v}_f$ is the Frobenius eigenvector of $\bld{A}^h$.  
\begin{figure}
\hspace*{\fill}
\begin{minipage}[b]{0.32\linewidth}
\centering
\includegraphics[width=\textwidth]{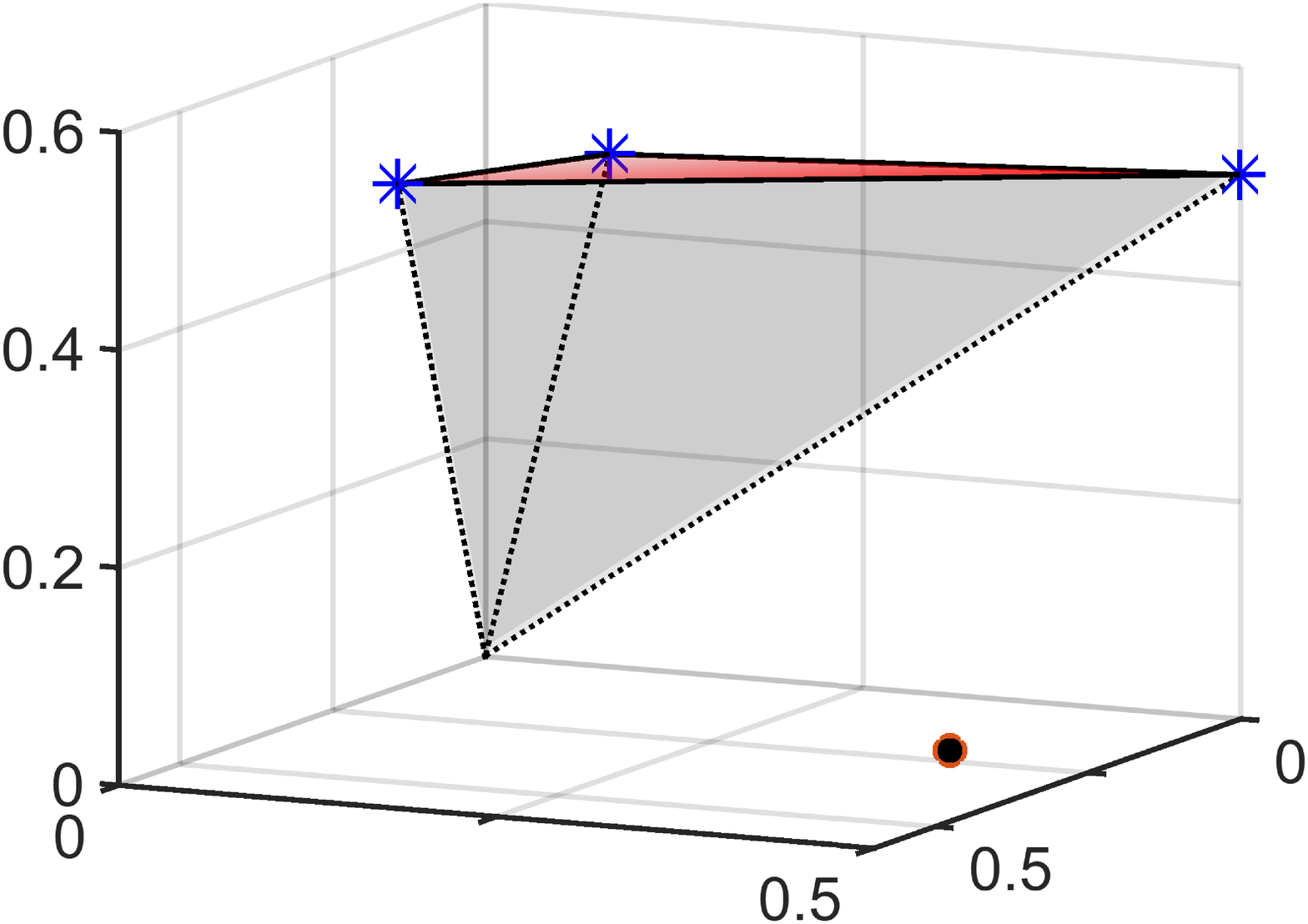}
\captionsetup{labelformat=empty}
\caption*{(a):~$\cset{3}{\bld{A},\bld{b}}$}
\end{minipage}
\hfill
\begin{minipage}[b]{0.32\linewidth}
\centering
\includegraphics[width=\textwidth]{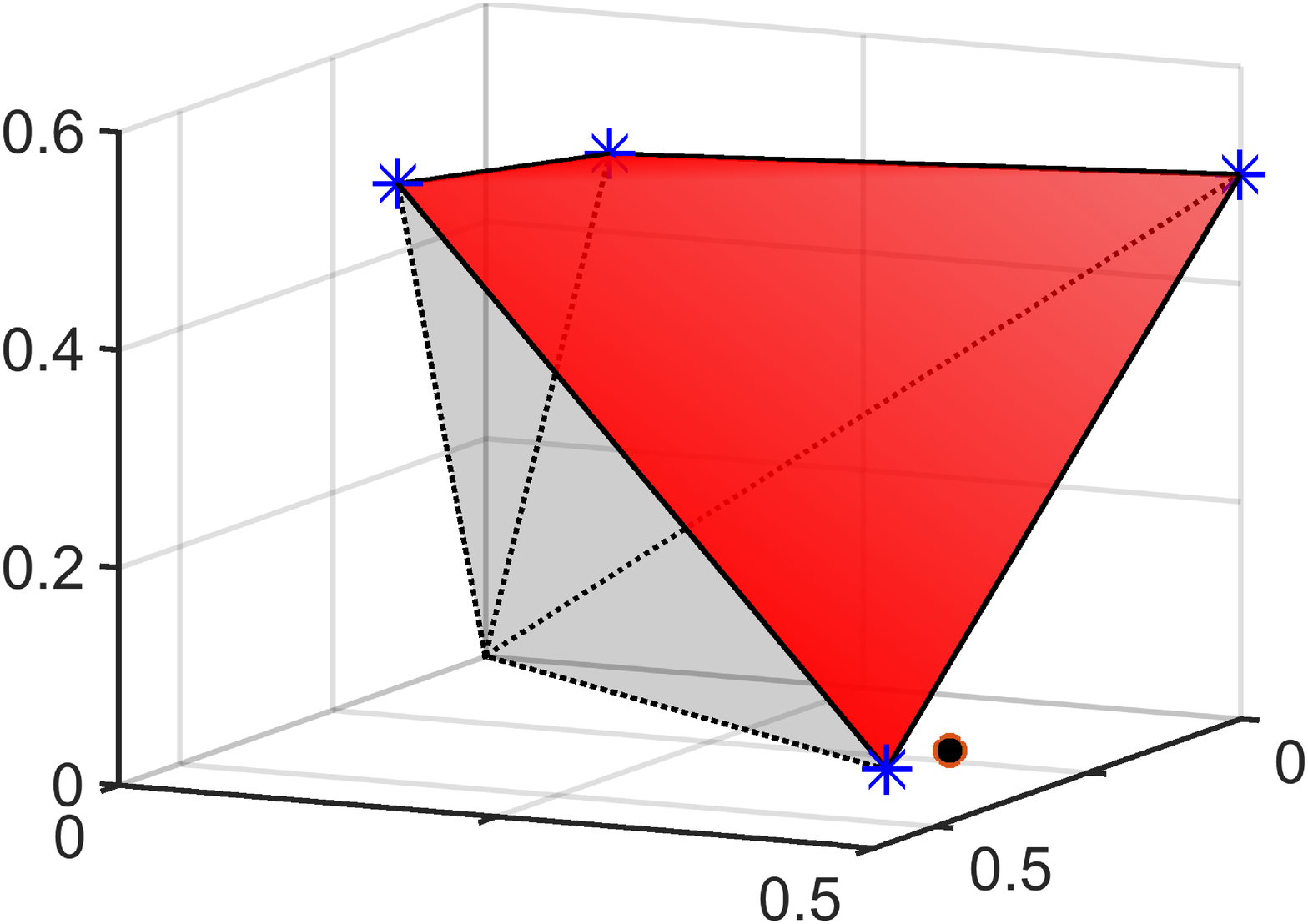}
\captionsetup{labelformat=empty}
\caption*{(b):~$\cset{4}{\bld{A},\bld{b}}$}
\end{minipage}
\hspace*{\fill}
\begin{minipage}[b]{0.32\linewidth}
\centering
\includegraphics[width=\textwidth]{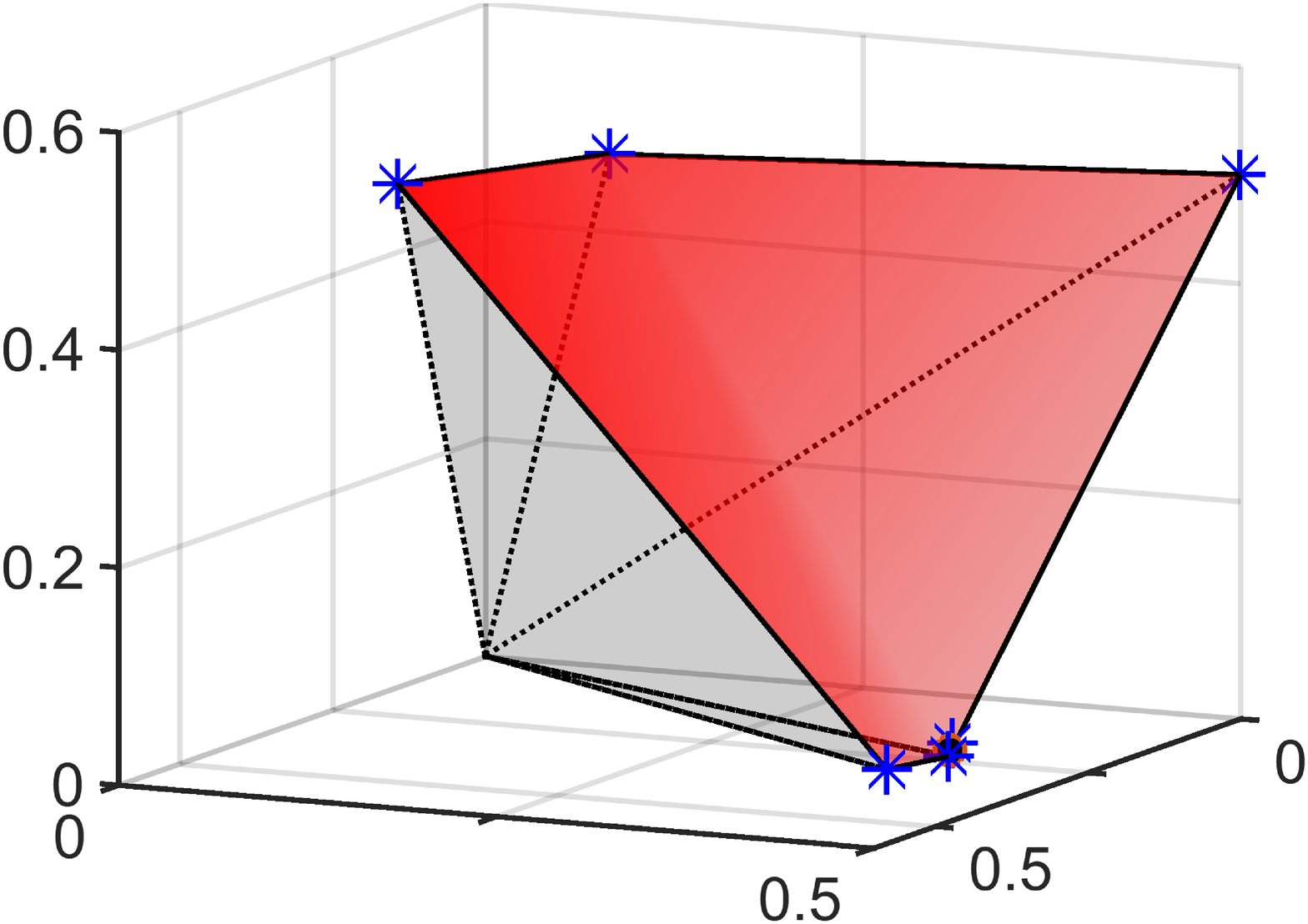}
\captionsetup{labelformat=empty}
\caption*{(c):~$\cset{6}{\bld{A},\bld{b}}$}
\end{minipage}
\hspace*{\fill}
\caption{\label{fig:E3}Example \ref{ex:3}: growth of the controllable cone mapped on the 3-dimensional simplex $S=\{\bld{x}\in \mathbb{R}^{3}_{+}|\trn{\mathbbm{1}}\bld{x}=1\}$; the generators of the cone and the Frobenius eigenvector are, respectively, marked by asterisks and a dot.}
\end{figure}
\end{exmp}
\subsection{Special Case}
So far it has been assumed that $\textup{rank}(\textup{conmat}_{n}(\bld{A},\bld{b}))=n$. Based on this assumption, the polyhedrality of the finite controllable subset only depends on the spectrum of $\bld{A}$. In addition, $k_{\textup{vert}}\geq n$ for $\cset{\textup{f}}{\bld{A},\bld{b}}$. We now point out that in the absence of such an assumption, $\cset{\textup{f}}{\bld{A},\bld{b}}$ can depend on the structure of $\bld{b}$ and that the vertex number can be less than $n$. In particular, it will be shown that $k_{\textup{vert}}=h$ if $\bld{b}\in\rn{n}$ is of a particular structure.
\begin{thm}
Let $\bld{A}\in \rn{n\times n}$ be irreducible with degree of cyclicity $h$ with\\ $0\leq h\leq n-1$. Then, $\cset{\textup{f}}{\bld{A},\bld{b}}=\textup{cone}\big{(}\textup{conmat}_{h}(\bld{A},\bld{b})\big{)}$ if \\ $\bld{b}\in\textup{cone}(\bld{v}_{f,0}\ \dots\ \bld{v}_{f,h-1})$, where $\bld{v}_{f,i}$, $i=0,\dots,h-1$ are the $h$ nonnegative eigenvectors of $\bld{A}^h$.
\end{thm}
\begin{pf}
Assume $\bld{b}=\displaystyle \sum_{i=0}^{h-1}c_{i}\bld{v}_{f,i}$ for some $\bld{c}\in \rn{h}$. Then, since \[\bld{A}^h\bld{b}=\displaystyle \sum_{i=0}^{h-1}c_{i}\rho(\bld{A})^h\bld{v}_{f,i}=\rho(\bld{A})^h\bld{b},\] it is immediate to see that $\bld{A}\big{(}\textup{conmat}_{h}(\bld{A},\bld{b})\big{)}=\big{(}\textup{conmat}_{h}(\bld{A},\bld{b})\big{)}\bld{X}$ has a nonnegative solution
\begin{equation}
\bld{X}=\begin{bmatrix}
0      & 0      & \cdots & 0 & \rho(\bld{A})^h\\
1      & 0      & \cdots & 0 & 0\\
0      & 1      &        & 0 & 0\\
\vdots &        & \ddots &   & \vdots\\
0      & \cdots & 0      & 1 & 0
\end{bmatrix},
\end{equation}
which, in the view of Corollary~\ref{cor:polyhedral_conset_f}, completes the proof.
\end{pf}
For $\bld{A}$ primitive (i.e., $h=1$), this results in the obvious case of $\cset{\textup{f}}{\bld{A},\bld{b}}$ being a ray along the Frobenius eigenvector $\bld{v}_{f}$ of $\bld{A}$ when $\bld{b}=c\bld{v}_{f}$ for any $c\geq 0$.
\section{Characterizations of Controllability}\label{sect:characterization_Controllability}
Given a cone $C_{\textup{obj}}\subseteq \rn{n}$ of control objectives or a subset of $\rn{n}$, the problem is to investigate whether $C_{\textup{obj}}$ is contained in $\cset{\textup{f}}{\bld{A},\bld{b}}$ or in $\cset{\infty}{\bld{A},\bld{b}}$. Of particular interest is when $C_{\textup{obj}}\subset \rn{n}$ is a polyhedral cone or a polytope. If the control objective cone $C_{\textup{obj}}$ is not polyhedral then outer approximate it by a polyhedral cone $C_{\textup{out}} \subseteq \mathbb{R}_+^n$ such that $C_{\textup{obj}} \subset C_{\textup{out}}$. Here, it is assumed that the controllabilty cone or its closure is polyhedral and that its corresponding vertex number or an upper bound of it is known. Hence $\cset{\infty}{\bld{A},\bld{b}}=\textup{cone}(\bld{b}\ \dots\ \bld{A}^{N-1}\bld{b}\ \bld{v}_{f,0}\ \dots\ \bld{v}_{f,h-1})$ for some $N\geq k_{\textup{vert}}^{\infty}$ and/or $\cset{\textup{f}}{\bld{A},\bld{b}}=\textup{cone}(\bld{b}\ \dots\ \bld{A}^{N-1}\bld{b})$ for some $N\geq k_{\textup{vert}}$. 
\begin{prop}\label{prop:Characterizaion_of_Controllability}
Let $C_{\textup{obj}}=\textup{conv}(\bld{p}_{1},\dots,\bld{p}_{m})$ or $C_{\textup{obj}}=\textup{cone}(\bld{p}_{1},\dots,\bld{p}_{m})$, where $\bld{p}_{i}\in \rn{n}$, $i=1,\dots,m$. 
\begin{enumerate}[(a)]
\item $C_{\textup{obj}}$ is controllable in finite time if and only if $\bld{p}\in \cset{\textup{f}}{\bld{A},\bld{b}}, \forall \bld{p}\in\{\bld{p}_{1},\dots,\bld{p}_{m}\}$. 
\item $C_{\textup{obj}}$ is controllable in infinite time (to be called {\em almost controllable}) if and only if $\bld{p}\in \cset{\infty}{\bld{A},\bld{b}}, \forall \bld{p}\in\{\bld{p}_{1},\dots,\bld{p}_{m}\}$ and $\exists\ \bld{p}^{\prime}\in \{\bld{p}_{1},\dots,\bld{p}_{m}\}$ such that $\bld{p}^{\prime}\notin \cset{\textup{f}}{\bld{A},\bld{b}}$. 
\end{enumerate}
\end{prop}
\begin{pf}
The proof is obvious from Definition~\ref{def:cset_f} and considering the fact that a cone can be expressed as a nonnegative combination of its generators.
\end{pf}
It is obvious from Proposition~\ref{prop:Characterizaion_of_Controllability}, that checking for controllability involves checking the following condition for each $i\in\{1,\dots,m\}$:
\begin{equation}\label{eq:Controllability_Solution}
\exists \bld{x}_{i}\in \{\bld{z}|\bld{Mz}=\bld{p}_{i},\bld{z}\in \rn{N}\},
\end{equation}
where $\bld{M}\in \rn{n\times N}$. Depending on the problem being investigated, either $\bld{M}=[\bld{b}\ \dots\ \bld{A}^{N-1}\bld{b}\ \bld{v}_{f,0}\ \dots\ \bld{v}_{f,h-1}]$ or $\bld{M}=[\bld{b}\ \dots\ \bld{A}^{N-1}\bld{b}]$.

In general, since $N\geq n$ (see Remark~\ref{rmk:polyhedral_cset_infty} and Remark~\ref{rmk:polyhedral_cset_f}), \Ref{eq:Controllability_Solution} defines an underdetermined system of equations. It is known that the nonnegative solution of  \Ref{eq:Controllability_Solution} is not unique in general \cite{Donoho_2005,Wang_2011}, and that uniqueness is guaranteed when the solution is sufficiently sparse \cite{Donoho_2005}. The authors of \cite{Donoho_2006} characterize necessary and sufficient conditions on the polytope $P=\textup{conv}(\bld{M})$ for uniqueness of the solution, where they prove unique solution exists if and only if $P$ is $k$-neighborly \footnote{A $k$-neighborly polytope is a convex polytope in which every set of $k$ or fewer vertices forms a face \cite{Convex_Polytops}.}. In \cite{Wang_2011,Donoho_2010}, the equivalent of this condition is presented in terms of the null space of $\bld{M}$. In this regard, this problem relates to the \emph{sparse measurement} problem, where it is formulated as reconstructing a nonnegative sparse vector from lower-dimensional linear measurements \cite{Khajehnejad_2011}. The results in this field do not directly apply here as the necessary sparsity condition is usually not met. In addition, we are not interested in finding the sparsest solution of \Ref{eq:Controllability_Solution}, which is normally an NP-hard problem \cite{Donoho_2005}.
\begin{prop}\label{prop:Characterization_of_Solution_Non-empty_Set}
Consider index sets $\mathcal{I}^{i}_{j}\subset \{1,\dots,N\}$ for $j=1,\dots,C(N,n)$ with $|\mathcal{I}^{i}_{j}|=n$, where $N>n$ is an upper bound to $k_{\textup{vert}}$ or an upper bound to $k_{\textup{vert}}^{\infty}$, $n$ is the dimension of space, and $C(N,n)$ is the number of $n$-combinations of the set $\{1,\dots,N\}$. Let $\bld{I}_{{\mathcal{I}}_{j}^{i}}$ denote the submatrix of the identity matrix of dimension $N$, $\bld{I}_N$ that is composed of the columns corresponding to ${\mathcal{I}}_{j}^{i}$. Then, \Ref{eq:Controllability_Solution} has a solution $\bld{x}_{i}$ for any $i\in \{1,\dots,m\}$ if and only if 
\begin{equation}\label{eq:Characterization_of_Solution_Non-empty_Set}
\bld{X}^{i}=\Big{\{}\bld{x}^{i}_{j}\big{|}\bld{x}^{i}_{j}=\bld{I}_{{\mathcal{I}}_{j}^{i}}(\bld{M}\bld{I}_{{\mathcal{I}}_{j}^{i}})^{-1}\bld{p}_{i},\ \bld{x}^{i}_{j}\in \rn{N},\ j=1,\dots,C(N,n)\Big{\}}
\end{equation}
is a non-empty set. 
\end{prop}
\begin{pf}
From our assumption we have $\bld{p}_{i}\in \textup{cone}(\bld{M})$. Since $N>n$, due to the Carath\'{e}odory theorem \cite{Caratheodory}, $\bld{p}_{i}$ also lies in at least one simplicial cone generated by $n$ columns of $\bld{M}$. Let $\mathcal{J}^i\in \{1,\dots,N\}$ with $|\mathcal{J}^i|=n$ be an index set composed of the indices of the columns generating this simplicial cone, and let $\bld{M}_{{\mathcal{J}}^i}$ denote the columns of $\bld{M}$ corresponding to $\mathcal{J}^i$. We can then write $\bld{p}_{i}\in \textup{cone}(\bld{M}_{{\mathcal{J}}^i})$, which can be expressed as $\bld{M}\bld{I}_{{\mathcal{J}}^i}\bld{z}^{i}=\bld{p}_{i}$ having a solution $\bld{z}^{i}\in \rn{n}$. Since $\bld{M}$ has full row rank and $\bld{I}_{{\mathcal{J}}^i}$ is full column rank, one obtains $\bld{z}^{i}=(\bld{M}\bld{I}_{{\mathcal{J}}^i})^{-1}\bld{p}_{i}$. Finally, we obtain a solution $\bld{x}_{j}^{i}\in \rn{N}$, where $\bld{x}_{j}^{i}=\bld{I}_{{\mathcal{J}}^i}\bld{z}^{i}=\bld{I}_{{\mathcal{J}}^i}(\bld{M}\bld{I}_{{\mathcal{J}}^i})^{-1}\bld{p}_{i}$. 

The converse is proved in a straightforward manner by noticing that every $\bld{z}\in \bld{X}^{i}$ satisfies \Ref{eq:Controllability_Solution}.
\end{pf}
\begin{rmk}\label{rmk:Controllability_Solution}
Let $\bld{X}^{i}=\Big{\{}\bld{x}^{i}_{1},\dots,\bld{x}^{i}_{q_{i}}\Big{\}}$ for some $q_{i}\in \mathbb{Z}_{+}$. It is then clear from the proof of Proposition~\ref{prop:Characterization_of_Solution_Non-empty_Set} that the set of solutions of \Ref{eq:Controllability_Solution} is the convex hull of $\bld{X}^{i}$, that is, we have for \Ref{eq:Controllability_Solution} that $\bld{x}_{i}\in \textrm{conv}(\bld{X}^{i})$.
\end{rmk}
Note that even though Proposition~\ref{prop:Characterization_of_Solution_Non-empty_Set} provides a method to determine whether $C_{\textup{obj}}\subseteq \textup{cone}(\bld{M})$ by checking inclusion of $C_{\textup{obj}}$ in any simplicial subcone of $\textup{cone}(\bld{M})$, the computational complexity of this method can be prohibitive as the check must be conducted for all $C(N,n)$ simplicial subcones in the worst case.
A more practical approach is then presented by the following proposition.
\begin{prop}\label{prop:LP_Solution}
Let $$\bld{M}_{f}=[\bld{b},\dots,\bld{A}^{N-1}b]\ \textup{and}\ \bld{M}_{\infty}=[\bld{b},\dots,\bld{A}^{N-1}b, \bld{v}_{f,0},\dots,\bld{v}_{f,h-1}].$$ Define the following optimization problem for each $i\in\{1,\dots,m\}$:
\begin{align}
\label{eq:Characterization_of_Controllability_LP}\underset{\bld{x}_{i}}{\min}&\ \trn{\bld{x}}_{i}\mathbbm{1}\\
\nonumber &\bld{M}\bld{x}_{i}= \bld{p}_{i}\\
\nonumber &\bld{x}_i\geq 0,
\end{align}
where $\mathbbm{1}\in \mathbb{R}^n$ is a vector of ones. We then have the following.
\begin{enumerate}[(a)]
\item The optimization problem \Ref{eq:Characterization_of_Controllability_LP} with $\bld{M}=\bld{M}_{\infty}$ has an optimal solution $\bld{x}^{\ast}_{i}\in \rn{N}$ if and only if \Ref{eq:Controllability_Solution} has a solution with $\bld{M}=\bld{M}_{\infty}$.
\item The optimization problem \Ref{eq:Characterization_of_Controllability_LP} with $\bld{M}=\bld{M}_{f}$ has an optimal solution $\bld{x}^{\ast}_{i}\in \rn{N}$ if and only if \Ref{eq:Controllability_Solution} has a solution with $\bld{M}=\bld{M}_{f}$. 
\end{enumerate}
\end{prop}
\begin{pf}
If \Ref{eq:Controllability_Solution} has a solution, the set $\bld{X}^{i}$ in \Ref{eq:Characterization_of_Solution_Non-empty_Set} is non-empty. As mentioned in Remark~\ref{rmk:Controllability_Solution}, the feasible set of \ref{eq:Characterization_of_Controllability_LP} is $\textup{conv}(\bld{X}^{i})$. Therefore, the convex optimization problem with linear penalty function converges to the minimum 1-norm solution in the feasible set. The converse is obvious.
\end{pf}
\begin{exmp}
We conclude this section with an example illustrating the application of Proposition~\ref{prop:LP_Solution}. Consider the system matrices of Example~\ref{ex:3}. Let $C_\textup{obj}$ be the polytope given by $$C_\textup{obj}=\Big{\{}\bld{p}\in \rn{4}\Big{|}\bld{p}=\sum_{i=1}^{4}\lambda_{i}\bld{p}_i, \lambda_{i}\geq 0, \sum_{i=1}^{4}\lambda_{i}=1\Big{\}},$$
where
\begin{align*}
\bld{p}_1  = [1,~3,~1,~1]^\textup{T},~&
\bld{p}_2  = [1,~3,~4,~3]^\textup{T},\\
\bld{p}_3  = [1,~2,~2,~1]^\textup{T},~&
\bld{p}_4  = [1,~1,~2,~1]^\textup{T}.                
\end{align*}
We will now check if the system initially at rest can be steered to any point in $C_{\textup{obj}}$ in finite time. From example~\ref{ex:3}, $k_{\textup{vert}}=6$ is known. Thus taking $\bld{M}=[\bld{b},~\bld{Ab},~\ldots,~\bld{A}^5\bld{b}]$, we solve for the optimization problem \Ref{eq:Characterization_of_Controllability_LP} using the Dual-Simplex algorithm implemented in Matlab Optimization Toolbox. The optimal solutions are obtained as
\begin{align*}
\bld{x}_1^{*}  = &~[0.1209,~0.3735,~0,~0.0078,~0,~0.0001]^\textup{T},\\
\bld{x}_2^{*}  = &~[2.3460,~0.6165,~0.0876,~0,~0.0003,~0]^\textup{T},\\
\bld{x}_3^{*} = &~[0.2989,~0.6982,~0.0473,~0,~0.0003,~0]^\textup{T},\\
\bld{x}_4^{*}  = &~[0.2517,~0.7798,~0.0071,~0,~0.0003,0]^\textup{T}.             
\end{align*}
Hence, the vertices of $C_{\textup{obj}}$ can be reached from the origin in finite number of steps using nonnegative inputs, which are determined by the solution vectors $\bld{x}_i^{*}$. Moreover, since $k_{\textup{vert}}=6$, every vertex of $C_{\textup{obj}}$ can be reached in at most 6 steps from the origin. Since $C_{\textup{obj}}$ is the convex hull of its vertices, we can conclude that any point $\bld{p}=\sum_{i=1}^{4}\lambda_{i}\bld{p}_i \in C_{\textup{obj}}$ can be reached from the origin in at most 6 steps using the input sequence $\bld{u}^{*}=\sum_{i=1}^{4}\lambda_{i}\bld{x}_i^{*}$.
\end{exmp}
\section{Conclusion}
We discussed a new view of the controllability problem for linear time-invariant positive systems that is more interesting for practical applications than the classical view. The controllabilty was defined as the ability to drive the system initially at origin to a certain target subset of $\rn{n}$ using nonnegative inputs. To this end, we discussed the geometry of controllable subsets and developed sufficient and necessary conditions for polyhedrality of such subsets. We showed that when the controllability matrix of the system is of full rank, those conditions solely depend on the spectrum of $\bld{A}$. In addition, it was shown that the controllable subset may keep growing for more than $n$ steps, where $n$ is the dimension of the system. We then proposed a numerical method to check for controllability of a linear positive system with respect to a certain objective set.

In this paper, we have focused on the single input case, where $\bld{b}\in \rn{n}$. The controllability problem for the multi-input case is an interesting problem as the results developed here are not directly applicable. The main issue, as noted in \cite{Benvenuti2006}, is that the direct sum of two non-polyhedral cones may result in a polyhedral cone. Therefore, one cannot apply the results of this paper to a set of system $(\bld{A},\bld{b}_i)$ separately, with $\bld{b}_i$ being a column of $\bld{B}$.\\
It is also of interest to investigate the geometry of controllable subsets when the controllability matrix is not of full rank. As far as the authors of this paper know, this is still an open issue.
\appendix
\newpage
\section{Positive Matrices}\label{sect:Appendix}
For completeness, we report Theorem~5 of \cite{Roitman1992} here. In this theorem,  $Q$ denotes the set of all real polynomials of the form $c_n x^n-\sum_{i=0}^{n-1}c_i x^i$, where $n\geq 1$, $c_n > 0$, and $c_i \geq 0$ for all $i$.
\begin{thm}[{\cite[Th.~5]{Roitman1992}}]\label{thm:Roitman}
Let $\{a_1,\dots,a_k\}$ be given complex numbers, and let $P(x)$ be the polynomial $x^{k} - a_{1}x^{k-l} - \dots - a_{k}$. Then conditions \Ref{item:A}, \Ref{item:B} and \Ref{item:C} below are equivalent:
\begin{enumerate}[(A)]
\item \label{item:A}Any infinite sequence $(u_{n})_{n\geq 0}$ of complex numbers which satisfies the recursion $u_{n+k}=a_{1}u_{n+k-1}+a_{2}u_{n+k-2}+\dots+a_{k}u_{n}$ for $n\geq 0$, also satisfies a recursion with nonnegative coefficients.
\item \label{item:B}The polynomial $P(x)$ divides a polynomial in $Q$.
\item \label{item:C}In case the polynomial $P(x)$ has a positive root $r$, then all conditions (1)-(4) below are satisfied:
\begin{enumerate}[(C1)]
\item \label{item:C1}$r\geq |\alpha|$ for any root $\alpha$ of $P(x)$.
\item \label{item:C2}if $\alpha=r$ for some root $\alpha$ of $P(x)$, then $\alpha/r$ is a root of unity.
\item \label{item:C3}all roots $P(x)$ with absolute value $r$ are simple.
\item \label{item:C4}if $P(r)=P(r\epsilon)=0$, where $\epsilon^k=1$ with $k\geq 1$ minimal, then $P(x)$ has no roots of the form $s\omega$ where $0<s<r$ and $\omega^k=~1$.
\end{enumerate}
\end{enumerate}
\end{thm}
\begin{lem}\label{lem:vf_existence}
Let $\bld{A}\in\rn{n\times n}$ be irreducible with cyclicity index $h$ and let $\bld{b}\in\rn{m}$. Define $C_{\lim}=\textup{cone}(\bld{A}_{f,0}\bld{b}, \ldots, \bld{A}_{f,h-1}\bld{b})$, where $\bld{A}_{f,i}=\displaystyle\textup{lim}_{k\rightarrow\infty}\frac{\bld{A}^{kh}}{\rho(\bld{A})^{kh}}\bld{A}^i$, for $i=0,\dots,h-1$. Let the nonnegative vectors $\bld{v}_{f,i}$, $i=0,\dots,h-1$ of Proposition~\ref{prop:polyhedral_conset_infty} be the $h$ distinct nonnegative eigenvectors of $\bld{A}^h$ associated with the Perron root of $\rho(\bld{A})^h$. It then holds that $C_{\lim}\subseteq\textup{cone}(\bld{v}_{f,0}\ \dots\ \bld{v}_{f,h-1})$, 
\end{lem}
\begin{pf}
Since $\bld{A}$ is irreducible, there exists a monomial matrix $\bld{S}\in \rn{n\times n}$ \cite[Th.~2.2.33]{Plemmons_book} , such that
\begin{equation}
\hat{\bld{A}}=\trn{\bld{S}}\bld{A}\bld{S}=\begin{bmatrix}
0_{n_1} &  A_{1} & 0 & \dots  & 0\\
\vdots & \ddots & \ddots & \ddots &\vdots\\
\vdots&  & \ddots  & \ddots & 0\\
0& \dots& \dots & 0_{n_{h-1}} & A_{h-1}\\
A_{h} & 0 &\dots &\dots & 0_{n_h}
\end{bmatrix},\
\hat{\bld{b}}=\trn{\bld{S}}b
\end{equation}
where $0_{n_i}\in \mathbb{R}^{n_i\times n_i}$, $i\in \mathbb{N}$ are square blocks with $\displaystyle\sum_{i=1}^{h}n_{i}=n$, and where $A_{i}$ has no zero rows or columns with $L_{1}=\displaystyle\prod_{i=1}^{h}A_{i}$ being an irreducible matrix. Then we have $\bld{\hat{A}}^h=\textup{diag}(L_{1},\dots,L_{h})$, where $L_{k}=\displaystyle\prod_{i=k}^h A_{i}\hspace{-4mm}\displaystyle\prod_{j=1}^{\textrm{mod}(h+k-1,h)}\hspace{-6mm}A_{j}$ is a primitive matrix of dimension $n_{k}\times n_{k}$ with Perron root $\rho(\bld{A})^h$. Define the matrix $\hat{\bld{A}}_{f,i}=\displaystyle\textup{lim}_{p\rightarrow\infty}\frac{\bld{\hat{A}}^{ph}}{\rho^{ph}}\bld{\hat{A}}^{i}$ for $i=0,\dots,h-1$. Since $L_{i}$, $i=1,\dots,h$ is primitive, it follows from \cite[Th.~2.4.1]{Plemmons_book} that 
\begin{equation}
\bld{\hat{A}}_{f,0}=\begin{bmatrix}
\bld{x}_{1}^{1}&\dots      &\bld{x}_{1}^{n_{1}}  &0           & 0       &0           &0   &\dots  & 0\\
0          &\dots      &0            &\bld{x}_{2}^{1} & \dots   &\bld{x}_{2}^{n_{2}} &0    &\dots  & 0\\ 
0          &\dots      &0            &0 & \dots   &0 &0    &\dots  & 0\\ 
\vdots     &\vdots     &\vdots       &\vdots      &\vdots   &\vdots      &\vdots&\vdots &\vdots\\
0          &0          &0            &0           &0        &0           &\bld{x}_{h}^{1}&\dots  &\bld{x}_{h}^{n_{h}}\end{bmatrix},
\end{equation}
where $\bld{x}_{i}^{k}=c_{i}^{k}\bld{x}_{i}$ with $c_{i}^{k}$, $k=1,\dots,n_{i}$, being some nonnegative scalars and with $\bld{x}_{i}\in \mathbb{R}_{s+}^{n_{i}\times n_{i}}$ being the Frobenius eigenvector of $L_{i}$. Note that due to the block structure of $\bld{\hat{A}}$, $\bld{\hat{A}}_{f,i}$ retains the same structure as $\bld{\hat{A}}_{f,0}$ up to a scaled permutation of its columns for $i=1,\dots,h-1$. Hence, we have $\bld{\hat{A}}_{f,i}\bld{\hat{b}}\in\textup{cone}(\bld{C})$, where
$$
\bld{C}=\begin{bmatrix}
\bld{x}_1 & 0         & \dots & 0\\
     0    & \bld{x}_2 & \dots & 0\\ 
     0    & 0         & \dots & 0\\
  \vdots  & \vdots    & \dots &\vdots\\
     0    & \dots     & 0     & \bld{x}_h
\end{bmatrix}.
$$
In the original coordinates, we have $\bld{A}_{f,i}\bld{b}\in\textup{cone}(\bld{SC})$. Clearly, since the columns of $\bld{C}$ are the nonnegative eigenvectors of $\bld{\hat{A}}^h$ and since $\bld{S}$ is monomial, we have $\bld{SC}=[\bld{v}_{f,0} \ \dots \ \bld{v}_{f,h-1}]$, where $\bld{v}_{f,i}\in \rn{n\times n}$ is the $(i+1)$-th nonnegative eigenvector of $\bld{A}^h$ for $i=0,\dots,h-1$. This proves that $\textup{cone}(\bld{A}_{f,0}b\ \dots\ \bld{A}_{f,h-1}b)\subseteq\textup{cone}(\bld{v}_{f,0}\ \dots\ \bld{v}_{f,h-1})$. 
\end{pf}
\begin{lem}\label{lem:subdominant_eigenvalues}
Let $\bld{A}\in\rn{n\times n}$ be irreducible with degree of cyclicity $h$ with $1\leq h\leq n$. Let $\bld{A}$ be decomposed as $\bld{A}=\bld{S}\textup{diag}(\bld{A}_1,\bld{A}_2)\bld{S}^{-1}$, where $\textup{spec}(\bld{A}_1)=\sigma^{\rho}(\bld{A})$ and $\textup{spec}(\bld{A}_2)=\sigma^{-}(\bld{A})$. Let $\sigma^{0}\subseteq \sigma^{\rho}(\bld{A}_2)$ be the set of all eigenvalues of $\bld{A}_2$ whose modulus is $\rho(\bld{A}_2)$ and whose polar angle is a rational multiple of $2\pi$. Then, there exists a minimal $M\in \mathbb{Z}_{+}$ such that 
\begin{equation}\label{eq:Appendix_rational_eigenvalues}
\sigma_0\subseteq \Big{\{}\lambda\in \textup{spec}(\bld{A}_2)\Big{|}\lambda=\rho(\bld{A}_2)\exp(\frac{2\pi k}{Mh} i), k=0,\ldots,Mh-1\Big{\}},
\end{equation}
or, equivalently, there exists a minimal $M\in \mathbb{Z}_{+}$ such that the eigenvalues of $\bld{A}_2/\rho(\bld{A}_2)$ with unit modulus whose argument are a rational multiple of $2\pi$ are among the $Mh$-th roots of unity.
\end{lem}
\begin{pf}
Let $\delta^{0}$ be a set of $n_{\delta^{0}}\in \mathbb{Z}_{+}$ members of $\sigma^{0}$ with the property that the difference between the polar angle of no two members of $\delta^{0}$ is an integer multiple of $2\pi/h$, or formally we define $\delta^{0}=\{\lambda_{1},\ldots,\lambda_{n_{\delta^{0}}}\in \sigma^{0}|\textup{arg}(\lambda_{i})-\textup{arg}(\lambda_{j})\neq 2z\pi/h, i\neq j, z\in \mathbb{Z}\}$. For $\lambda_{j}\in \delta^{0}$, $j=1,\ldots,n_{\delta^{0}}$, let $\textup{arg}(\lambda_{j})=\frac{2\pi p_j}{q_j}$. Define the sets $\sigma^{0}_{j}\subset \sigma^{0}$ for $j=1,\ldots,n_{\delta^{0}}$ as
$$\sigma^{0}_{j}=\Big{\{}\lambda\in \textup{spec}(\bld{A}_2)\Big{|}\lambda=\rho(\bld{A}_2)\exp\big{(}(k/h+p_j/q_{j})2\pi i\big{)},~k=0,\ldots,h-1\Big{\}},$$
or equivalently using the notation $s_{j,k}\equiv kq_j+hp_{j} (\textup{mod}~hq_j)$,
$$\sigma^{0}_{j}=\Big{\{}\lambda\in \textup{spec}(\bld{A}_2)\Big{|}\lambda=\rho(\bld{A}_2)\exp\big{(}\frac{s_{j,k}}{hq_{j}}2\pi i\big{)},~k=0,\ldots,h-1\Big{\}}.$$
It is clear that $\sigma^{0}_{1},\ldots,\sigma^{0}_{n_{\delta^{0}}}$ are mutually disjoint. In addition, since the eigenvalues of $\bld{A}$ are invariant under polar rotation of $2k\pi/h$ for any $k\in \mathbb{Z}$, we have $\sigma^{0}=\cup_{j=1}^{n_{\delta^{0}}}\sigma^{0}_{j}$. Noting that $0\leq s_{j,k}\leq hq_{j}-1$ for $k=0,\ldots,h-1$ and for $j=1,\ldots,n_{\delta^{0}}$, one observes that $\sigma_0$ has the from proposed in \Ref{eq:Appendix_rational_eigenvalues} by choosing $M=\textup{lcm}(q_1,\ldots,q_{n_{\delta^{0}}})$.
\end{pf}
\newpage
\section*{References}
\bibliography{PhDref_ver2.bib}
\end{document}